\documentclass{amsart}
\usepackage{amscd}
\usepackage{amssymb}
\usepackage[dvips]{epsfig}

\newcommand{\al}{\alpha}
\newcommand{\Ap}{\Delta_{K}(t)}
\newcommand{\Apn}{\Delta_{K_n}(t)}
\newcommand{\Apone}{\Delta_{K_n}(-1)}
\newcommand{\B}{\mathcal B}
\newcommand{\brho}{\bar{\rho}}
\newcommand{\brhoz}{\bar{\rho}_0}
\newcommand{\bX}{\bar{X}}
\newcommand{\C}{\mathbb C}
\newcommand{\card}[1]{\mbox{card}( #1 )}
\newcommand{\coker}[1]{\mbox{coker} ( #1 )}
\newcommand{\CP}{\mathbb{CP}^1}
\newcommand{\dmn}[1]{\mbox{dim}( #1 )}
\newcommand{\dimC}[1]{\mbox{dim}_{\C}( #1 )}
\newcommand{\dn}{\Delta(2,3,|n|)}
\newcommand{\dpqr}{\Delta(p,q,r)}
\newcommand{\dpt}{\pi_1(\partial \Mt)}
\newcommand{\fl}[1]{\lfloor #1 \rfloor}\newcommand{\gl}{\mathfrak{l}}
\newcommand{\fltwo}[1]{\lfloor \frac{#1}{2} \rfloor}
\newcommand{\hB}{\widehat{B}}
\newcommand{\hM}{\widehat{M}}
\newcommand{\hp}{\widehat{\pi}}
\newcommand{\hS}{\widehat{\Sigma}} 
\newcommand{\I}{I}
\newcommand{\img}[1]{\mbox{Im} ( #1 )}
\newcommand{\lt}{\tilde{\lambda}}
\newcommand{\m}{\mathfrak{m}}
\newcommand{\mat}[4]{\left( \begin{array}{cc} #1 & #2 \\
    #3 & #4 \end{array} \right)}
\newcommand{\Md}{\Mt_d}
\newcommand{\mt}{\tilde{\mu}}
\newcommand{\Mt}{\widetilde{M}}
\newcommand{\noy}[1]{\mbox{ker} ( #1 )}
\newcommand{\pB}{\pi_1^{\mbox{orb}}(\B)}
\newcommand{\pdM}{\pi_1(\pM)} 
\newcommand{\pM}{\partial M}
\newcommand{\Pmat}[4]{ \pm \left( \begin{array}{cc} #1 & #2 \\
    #3 & #4 \end{array} \right)}
\newcommand{\pn}{\pi}
\newcommand{\pd}{\pi_1(\St)}
\newcommand{\PSLC}{\mbox{PSL}_2(\C)}
\newcommand{\pt}{\tilde{\pi}}
\newcommand{\Q}{\mathbb Q}
\newcommand{\R}{\mathbb R}
\newcommand{\rhod}{\rho_d}
\newcommand{\rhoz}{\rho_0}
\newcommand{\Sd}{\Sigma_d}
\newcommand{\slC}{sl_2(\C)}
\newcommand{\SLC}{\mbox{SL}_2(\C)}
\newcommand{\SLR}{\mbox{SL}_2(\R)}
\newcommand{\St}{\Sigma_2}
\newcommand{\tr}[1]{\mbox{trace}( #1 )}
\newcommand{\tXi}{\widetilde{X}_i}
\newcommand{\Xin}{X_i^{\nu}}
\newcommand{\Z}{\mathbb Z}

\newtheorem{theorem}{Theorem}[section]
\newtheorem{claim}[theorem]{Claim}
\newtheorem{con}[theorem]{Conjecture}
\newtheorem{defn}[theorem]{Definition}
\newtheorem{lemma}[theorem]{Lemma}
\newtheorem{prop}[theorem]{Proposition}
\newtheorem{que}[theorem]{Question}
\newtheorem{sub}[theorem]{Subclaim}
\newcommand{\Pf}{\noindent {\bf Proof: }}
\newcommand{\Rmk}{\noindent {\bf Remark: }}

\begin{document}
\title{The Culler-Shalen seminorms of the $(-2,3,n)$ pretzel knot}
\date{}
%\date{August 8, 2001}

\author{Thomas W. Mattman}
\address{Department of Mathematics \& Statistics, 
California State University, Chico\\
Chico, CA95929-0525}
\email{TMattman@CSUChico.edu}
\thanks{Research supported in part by grants from NSERC (Canada),
FCAR (Qu\'ebec), and RIMS (Kyoto)}

\begin{abstract}
We show that the $\SLC$-character variety of the $(-2,3,n)$ pretzel 
knot consists of two (respectively three) algebraic curves 
when $3 \nmid n$ (respectively $3 \mid n$) and give an
explicit calculation of the Culler-Shalen seminorms 
of these curves. Using this calculation, we describe
the fundamental polygon and Newton polygon for these knots 
and give a list of Dehn surgeries yielding a manifold with
finite or cyclic fundamental group. This constitutes a new proof
of property P for these knots.
\end{abstract}

\maketitle

\section*{Introduction}

Let $M = S^3 \setminus K$ denote the exterior of a hyperbolic knot $K$.
In \cite{CGLS}, \cite{CS1}, and \cite{CS2}, Culler and Shalen construct a 
norm on the vector space $V = H_1(\partial M ;{\R})$ 
which is a powerful
tool in the study of Dehn surgery on $K$. In particular, they show 
that if surgery along slope $\alpha$ results in a manifold $M(\alpha)$ 
having cyclic fundamental group and $\alpha$ is not a boundary slope, 
then $\alpha$ has minimal norm. 
Boyer and Zhang~\cite{BZ} extended this work by arguing that if
$\alpha$ is not a boundary slope and $\pi_1(M(\alpha))$ is finite, then, 
again, $\alpha$ must have small norm. 
Taking advantage of this observation, we~\cite{BMZ} proved that 
the $(-2,3,7)$ pretzel knot admits only four finite
surgeries and that its $\SLC$-character variety consists of two algebraic
curves. In the current article we   
generalize these results to hyperbolic pretzel knots of the form $(-2,3,n)$. 

In particular, we argue that 
the only examples of non-trivial 
cyclic or finite surgeries on such a 
knot are the five surgeries on 
the $(-2,3,7)$ and $(-2,3,9)$ 
pretzels found by Fintushel and Stern~\cite{FS} 
and Bleiler and Hodgson~\cite{BH}.
Thus, we provide an alternate proof of Delman's~\cite{D} result that 
non-trivial surgeries on these knots yield 
manifolds with infinite fundamental group when $n < 0$. 
Moreover, we extend this observation to the case $n \geq 11$.
As for the number of curves in the $\SLC$-character variety,
we prove that there are two when $3 \nmid n$ and three otherwise.

The essential new ingredient is the use of the Seifert fibred 
surgeries of the $(-2,3,n)$ pretzel knots. Although the 
Culler-Shalen norm of 
a slope $\alpha$ may be very large when $M(\alpha)$ is a 
Seifert fibred space, it is nonetheless possible to bound
that norm in terms of the Seifert indices of $M(\alpha)$
(see \cite{BB}). Similarly, the Seifert structure of the 
$2$-fold branched cyclic cover of a pretzel knot also constrains the
norm. Combining these constraints with knowledge of the 
boundary slopes (\cite{HO,Du}) allows us to explicitly work out
the Culler-Shalen norms of the hyperbolic $(-2,3,n)$ pretzel knots.
Given these norms, we classify the 
finite and cyclic Dehn surgeries, enumerate the  
curves in the $\SLC$-character variety,
and construct the fundamental polygon and Newton polygon of the
$A$-polynomial for these knots.

Section 1 provides
some preliminary definitions and a brief review of the 
theory of Culler-Shalen seminorms.
With these in hand, we state our results more explicitly and
conclude that section with an outline of the paper.

\section{Preliminaries, Results, and Outline}

\subsection*{The $(-2,3,n)$ Pretzel Knot}

\begin{figure}
\begin{center}
\epsfig{file=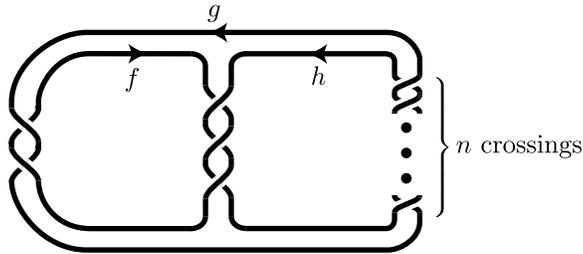}
\end{center}
\caption{The $(-2,3,n)$ Pretzel Knot \label{fgprtz}}
\end{figure}

Let $K_n$ denote the $(-2,3,n)$ pretzel knot (see Figure~\ref{fgprtz}).
If $n$ is even, $K_n$ is a link.
Also, $K_1$, $K_3$ and $K_5$ are torus knots and therefore not hyperbolic.
So we will assume that $n$ is an odd integer, $n \neq 1,3,5$. 

Let $\pn$ denote the fundamental group of $M = S^3 \setminus K_n$
and $\pt$ that of its $2$-fold cover $\Mt$. The $2$-fold branched
cyclic cover will be denoted by $\St$ and we will make strong use of
the fact that $\pd = \pt /  \langle \mu^2 \rangle = 
\pt / \langle \tilde{\mu} \rangle$ where
$\mu \in \pi$ is the class of a meridian of $K_n$ 
and $\tilde{\mu} \in \pt$ the class of the loop in $\Mt$ 
which (double) covers that meridian. 
Similarly $\lambda$ will denote the class of a preferred longitude of $K_n$ 
and $\lt \in \pt$ its lift. At the same time,
$\pd$ is a central extension of the triangle group $\dn$.
This is because $\St$ is Seifert fibred with base orbifold
$\B = S^2(2,3,|n|)$ and $\pB = \dn$.

The manifold $M$ is {\it small} in the sense that it contains
no closed essential surfaces~\cite{O}.  Essential
surfaces therefore meet $\pM$ in a non-empty set of parallel curves
each having the same slope. Slopes which can be obtained in this manner
are called {\it boundary} slopes.

\subsection*{The Character Variety and Culler-Shalen Seminorms}

We refer the reader to \cite[Chapter 1]{CGLS} and \cite{BZ2} for a more
detailed exposition. In particular, we restrict 
ourselves here to the case of small hyperbolic knots.

Let $R = \mbox{Hom}(\pn, \SLC)$ denote the set of
$\SLC$-representations of the fundamental group of $M$.
Then $R$ is an affine algebraic set, as is $X$, the set of characters
of representations in $R$.
Since $M$ is small, the irreducible components of $X$ are 
curves~\cite[Proposition 2.4]{CCGLS}. Moreover, for each component
$R_i$ of $R$ which contains an irreducible representation, the 
corresponding curve $X_i$ induces a non-zero seminorm $\| \cdot \|_i$ on 
$V = H_1(\pM; \R)$~\cite[Propositon 5.7]{BZ2} via the following construction.

For $\gamma \in \pn$, define the regular function
$I_{\gamma}:X \to \C$ by $I_{\gamma}(\chi_{\rho}) = \chi_{\rho}(\gamma) =
\mbox{trace}(\rho(\gamma))$.
By the Hurewicz isomorphism,
a class $\gamma \in L =  H_1(\pM ;{\Z})$ determines an element
of $\pi_1(\pM)$, and therefore an element of $\pn$ well-defined
up to conjugacy. The function
$f_{\gamma} = I_{\gamma}^2 -4$ is again regular and so can be pulled back to
$\tXi$, the smooth projective variety birationally equivalent to
$X_i$. For $\gamma \in L$, $\| \gamma \|_i$ is the degree of
$f_{\gamma} : \tXi \to \CP.$ The seminorm is extended to $V$ 
by linearity. We will call a seminorm constructed in this manner 
a Culler-Shalen seminorm.

If no $f_{\gamma}$ is constant on $X_i$, then $\| \cdot \|_i$ is in fact a
norm (rather than just a seminorm) 
and we shall refer to $X_i$ as a {\em norm curve}. 
Since $M$ is hyperbolic, there is a representation into $\PSLC$ carrying 
the hyperbolic structure. This lifts to an $\SLC$-representation
and the curve on which the character of this representation lies
is a norm curve \cite{CGLS}.
We will refer to it as the {\em canonical curve}.
If $X_i$ is not a norm curve, then there is a boundary slope 
$r$ such that $f_r$ is constant on $X_i$.
In this case, we will call $X_i$ an {\em $r$-curve}.

The minimal norm 
$s_i = \min \{ \| \gamma \|_i \, ; \, \gamma \in L, \, \| \gamma \|_i > 0 \}$ 
is an even integer, as is
$S = \sum_{i} s_i$, the sum being taken over the curves $X_i \subset X$. 
We will denote the sum of the Culler-Shalen 
seminorms by $\| v \|_T$ (here $v \in
V$), i.e.,  $\| v \|_T = \sum_i \| v  \|_i$.
Since the sum includes the norm on the canonical curve, $\| \cdot \|_T$ 
is a norm (not just a seminorm).

The power of Culler-Shalen seminorms is illustrated by the following 
two theorems which relate them to finite and cyclic surgeries.

\setcounter{section}{0}

\begin{theorem}[Corollary 1.1.4 \cite{CGLS}] \label{Thcyc}
If $\alpha$ is not a boundary slope and $\pi_1(M(\alpha))$ is cyclic, 
then $\| \alpha \|_i = s_i$. 
\end{theorem}

\begin{theorem}[Theorem 2.3 \cite{BZ}] \label{Thfin}
If $\alpha$ is not a boundary slope and $\pi_1(M(\alpha))$ is finite, 
then $\| \alpha \|_i \leq \mbox{max}(2s_i, s_i + 8)$. 
\end{theorem}

An important example is {\em trivial} surgery, i.e., surgery along
the meridian $\mu$. 
Since $\mu$ is not a boundary slope of $K_n$ \cite{O} and 
$M(\mu) = S^3$ has cyclic fundamental group, 
$\| \mu \|_i = s_i$ for each 
Culler-Shalen seminorm and $\|\mu\|_T = S$.

We will also be making use of the strong connection with boundary slopes.
Using meridian-longitude coordinates, the slopes of $K_n$ are parameterized
by $\Q \cup \{ 1/0 \}$. The {\em distance} 
$\Delta(a/b, c/d)$ between two such slopes $a/b$ and $c/d$ is their
minimal geometric intersection number $|ad-bc|$. In the case of a 
knot, such as $K_n$, for which $\mu$ is not a boundary slope, we can
rewrite Lemma~6.2 of \cite{BZ} as follows. (Note that 
the underlying idea is implicit in \cite[Chapter 1]{CGLS}.)

\begin{lemma}[Lemma 6.2 \cite{BZ}] \label{lemBZ62}
\begin{equation*} 
\| \gamma \|_i = 2 {[} \sum_j a^i_j \Delta( \gamma, \beta_j) {]}
\end{equation*}
where the $a^i_j$ are non-negative integers and the 
sum is over the set of boundary slopes $\beta_j$.
\end{lemma}  

On a norm curve $X_i$, at least two of the $a^i_j$ are non-zero. 
On an $r$-curve, $r=\beta$ is a boundary slope, and only the 
$a^i_j$ corresponding to $\beta$ is non-zero. Thus, the 
Culler-Shalen seminorm
on an $r$-curve is of the form $\| \gamma \|_i = s_i \Delta( \gamma, r)$.
Since $\mu = 1/0$ must have minimal norm $s_i$, it follows that
$r$ is an integral boundary slope.

\setcounter{section}{1}
\setcounter{theorem}{0}

\subsection*{Results and Outline of the paper}

Our results rest on Propositions~\ref{prpext} and
\ref{prpS} both of which apply to more general classes 
of knots (with only minor changes to the proofs). Therefore, these may be of 
independent interest.

\begin{prop} \label{prpext} 
Let $\brhoz$ be an irreducible $\PSLC$-representation of
$\pt$ which factors through $\dn$. Then $\brhoz$ has a unique
extension to $\pi$.
\end{prop}

This can be generalized to three-tangle Montesinos knots 
$\m =  K(a/p, b/q, c/r)$ 
by replacing $\dn$ with the triangle group $\dpqr$
(see \cite{Ma1} for details). 

\medskip

\Rmk Note that $\dpqr$ is the orbifold fundamental group of 
$\B$, the base orbifold of the $2$-fold branched cyclic 
cover of $\m$. Essentially, Proposition~\ref{prpext}
says the $\PSLC$ character variety of $\pB = \dpqr$ includes into that of 
the knot $\m$: $\bX(\pB) \subset \bX(\m)$. 

Although we cannot prove such an inclusion for more general Montesinos
knots, this nonetheless suggests that the character variety of a
Montesinos knot $\m$ is largely determined by that of its associated
orbifold $\B$. In particular, note that 
the dimension of the $\SLR$ character variety $X_{\R}$ of $\pB$ grows linearly
with the number of tangles in $\m$. Indeed, as a real variety,
$X_{\R}$ includes the Teichm\"uller space of $\B$ which is 
isomorphic to $\R^{2(t-3)}$ where $t$ is the number of tangles in 
$\m$. Based on this evidence, we conjecture that the character varieties
of $\m$ show the same behaviour.

\begin{con} The dimensions of the $\PSLC$- and $\SLC$-character varieties
of Montesinos knots $\m = \m(a_1/b_1, a_2/b_2, \ldots, a_t/b_t)$ grow 
linearly with the number of tangles $t$.
\end{con}

\begin{prop} \label{prpS}
The minimum of the total norm $\| \cdot \|_T$ is $S = 3(|n-2| -1)$. 
\end{prop}
More generally, for a $(-2,p,q)$ pretzel knot (see \cite{Ma1}), 
$$S = |pq| - (|p| + |q|) + |pq - 2(p+q)|.$$ 

The specific property of the $(-2,3,n)$ pretzel knots used in 
our argument is the existence of two Seifert fibred surgeries
at slopes $2n+4$ and $2n+5$, first observed by Bleiler and Hodgson~\cite{BH}.
Using the work of Ben Abdelghani and Boyer~\cite{BB}, we can determine the 
norm of these slopes.

\begin{prop} \label{prp2n4}
The total norm of the $2n+4$ Seifert fibred surgery is
$\| 2n+4 \|_T = S + 3(|n-6| - 1)$.
\end{prop}

\begin{prop} \label{prp2n5}
The total norm of the $2n+5$ Seifert fibred surgery is
$\| 2n+5 \|_T = S + 4(|n-5| - 2)$.
\end{prop}

A final ingredient is the connection with boundary slopes
(Lemma~\ref{lemBZ62}). 
The boundary slopes can be found using the methods of 
\cite{HO,Du} and with those in hand, the calculation of
the seminorms comes down to determining the integers $a^i_j$
of Lemma~\ref{lemBZ62}.

A careful analysis of the possible values for the $a^i_j$ 
allows us to deduce our Main Theorem.

\begin{theorem} \label{thmM} The $\SLC$ character variety of the 
hyperbolic $(-2,3,n)$ pretzel knot $K_n$ contains a curve
of reducible characters and a norm curve $X_0$. If $n \mid 3$,
there is in addition an $r$-curve $X_1$ with $r = 2n+6$ and
$s_1 = 2$. The Culler-Shalen norm $\| \cdot \|_0$ 
for the norm curve is as follows.

If  $3 \nmid n$, then $s_0 = 3(|n-2| -1)$ and
$$ \|\gamma\|_0 = 2 [ \Delta (\gamma,16) + 
2 \Delta( \gamma, \frac{n^2 - n - 5}{ \frac{n-3}{2}}) 
+ \frac{n-5}{2} \Delta( \gamma, 2n+6)]$$
when $n \geq 7$ and 
$$ \|\gamma\|_0 = 2 [  \Delta (\gamma,10)  
+ \frac{1-n}{2} \Delta( \gamma, 2n+6) + 
\Delta( \gamma, 2(n+1)^2 / n) ]$$
when $n \leq -1$.
 
If $3 \mid n$, then 
$s_0 = 3|n-2| -5$ and 
$$ \|\gamma\|_0 = 2 [ \Delta (\gamma,16) + 
2 \Delta( \gamma, \frac{n^2 - n - 5}{ \frac{n-3}{2}}) 
+ \frac{n-7}{2} \Delta( \gamma, 2n+6)]$$
when $n \geq 7$ and 
$$ \|\gamma\|_0 = 2 [  \Delta (\gamma,10)  
- \frac{n+1}{2} \Delta( \gamma, 2n+6) + 
\Delta( \gamma, 2(n+1)^2 / n) ]$$
when $n \leq -1$.  
\end{theorem}

Once we have this theorem, we can determine which slopes of $K_n$
are candidate cyclic or finite surgery slopes as these are either 
of small norm, or else boundary slopes (Theorems~\ref{Thcyc} and \ref{Thfin}).
We conclude that 
the only non-trivial cyclic or finite surgeries on these $K_n$ 
are the
five on the $(-2,3,7)$ and $(-2,3,9)$ pretzel knots discovered
by Fintushel and Stern (see \cite{FS}) and Bleiler and Hodgson~\cite{BH}. 
(Indeed,
these are the only non-trivial finite or cyclic surgeries on
any non-torus $(p,q,r)$ pretzel knot, see~\cite{Ma2}.) More precisely,
we have

\begin{prop} \label{prp23n}
If the $(-2,3,n)$ pretzel knot $K_n$ admits a non-trivial
cyclic or finite surgery, then one of the following holds.
\begin{itemize}
\item $K_n$ is torus, in which case $n = 1$,$3$, or $5$,
\item $n =7$, in which case $18$ and $19$ are cyclic
fillings while $17$ is a finite, non-cyclic filling, or
\item $n=9$, in which case $22$ and $23$ are finite, non-cyclic
fillings.
\end{itemize}
\end{prop}

We can also use Theorem~\ref{thmM} to determine
the fundamental polygon $B$ which is the disc of radius $s_0$ in
$H_1(\pM ;R)$. As $B$ is dual to the Newton polygon $N$ of the 
$A$-polynomial \cite{CCGLS}, we can likewise describe $N$ explicitly. 
This is significant as it remains difficult to calculate 
$A$-polynomials of knots.

In summary then, there are three main inputs for our approach.
We use information about the two-fold branched cyclic cover and
the boundary slopes of our knots. For example, these are both
known for the Montesinos knots. We also take advantage of
the two Seifert fibred slopes $2n+4$ and $2n+5$. We parlay this
data into information about the $\SLC$-character variety and
about cyclic and finite surgeries of the knot.
In principle this procedure
could be carried out for any Montesinos knot which admits a
Seifert fibred (or finite or cyclic) surgery. Indeed, in
our thesis~\cite{Ma1} we make a similar analysis of the 
twist knots and the $(-3,3,n)$ pretzel knots $(|n| \leq 6)$
(see also \cite{BMZ}).
This leads us to ask 

\begin{que} Are there other examples of Montesinos knots
admitting Seifert fibred (or cyclic or finite) surgeries?
\end{que}
We expect that our method could be applied to such examples.

The structure of our paper is as follows. 
Propositions~\ref{prpext}, \ref{prpS}, \ref{prp2n4},
and \ref{prp2n5} are proved in Sections 2, 3, 4, and 5 respectively.
These propositions are then used to prove our Main 
Theorem~\ref{thmM} in Section 6. In Section 7 we
use the theorem to prove Proposition~\ref{prp23n} and to
describe the fundamental polygons and
Newton polygons of the hyperbolic $(-2,3,n)$ 
pretzel knots.

\section{Proof of Proposition~\ref{prpext}}

In this section we prove

\setcounter{section}{1}

\begin{prop}  
Let $\brhoz$ be an irreducible $\PSLC$-representation of
$\pt$ which factors through $\dn$. Then $\brhoz$ has a unique
extension to $\pi$.
\end{prop}

\setcounter{section}{2}
\setcounter{theorem}{0}

\noindent{\bf Notation:} In general $\rho$ will denote an $\SLC$-representation
of $\pi$ and $\rhoz$ its restriction to $\pt$. The corresponding
$\PSLC$-representations will be denoted $\brho$ and $\brhoz$ respectively.

\medskip

\Pf
Suppose $\brho$ and $\brho^{\prime}$ were two extensions. Let 
$\alpha \in \pi \setminus \pt$. For any $\beta \in \pt$, 
$A = \brho(\alpha)^{-1} \brho^{\prime}(\alpha)$ commutes with $\brhoz(\beta)$.
So $A$ commutes with each element of $\brhoz(\pt)$.
However, as $\brhoz$ is irreducible, this implies
$A = \pm I$. Thus, $\brho$ and $\brho^{\prime}$ 
agree on $\alpha$ and hence on $\pi$. 

Since $\pi = \pt \sqcup \mu \pt$, there will be an extension 
provided we can find $A \in \PSLC$ with $A^2 = \pm I$ and such that
$A \brhoz(\beta) A^{-1} = \brhoz(\mu \beta \mu^{-1})$ for all
$\beta \in \pt$. Our goal then is to find such an $A$ corresponding
to conjugation by $\mu$.

\begin{figure}
\begin{center}
\epsfig{file=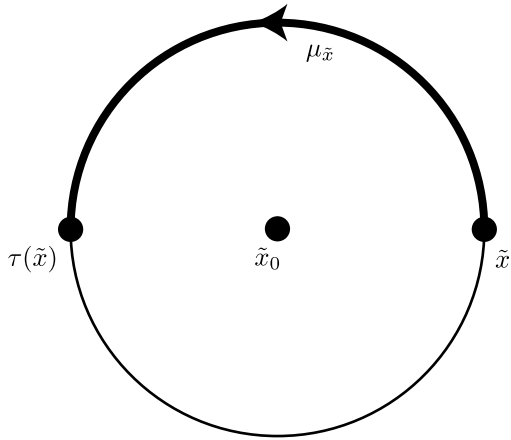}
\caption{\label{fgmux}%
A lift of the meridian $\mu$.}
\end{center}
\end{figure}

Let $\tau$ be the involution of the 2-fold branched cyclic
cover $\St$ and choose base points $\tilde{x} \in \partial \Mt$
and $\tilde{x}_0 \in \mbox{Fix}(\tau) \subset \St$.
Then conjugation by $\mu$ corresponds to 
$$ \begin{array}{ccc}
\pi_1(\Mt, \tilde{x}) & \stackrel{\cong}{\longrightarrow} & 
\pi_1(\Mt, \tilde{x}) \\
{[} \alpha {]}   & \mapsto & [ \alpha ]^{\mu} = 
[ \mu_{\tilde{x}} \cdot \tau(\alpha) \cdot (\mu_{\tilde{x}})^{-1}] 
\end{array} $$
where $\mu_{\tilde{x}}$ is the lift of $\mu$ beginning at $\tilde{x}$
(see Figure~\ref{fgmux}).

\begin{figure}
\begin{center}
\epsfig{file=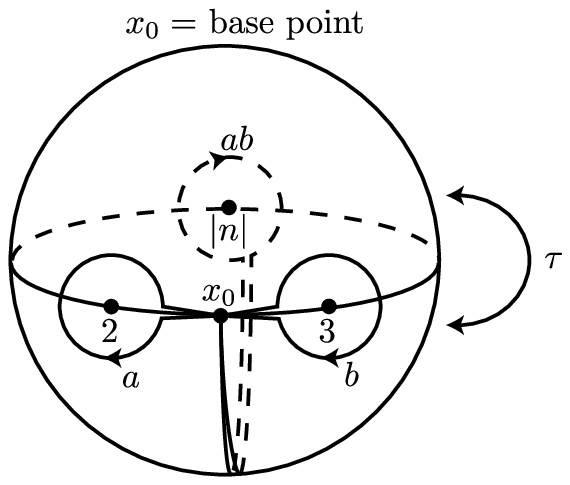}
\end{center}
\caption{$\tau$ action on $S^2(2,3,|n|)$ \label{fgtact}}
\end{figure}

On the other hand, restricted to the base orbifold 
$\B = S^2(2,3,|n|)$, $\tau$ is reflection in the equator 
and interchanges the hemispheres
(see~\cite[Th\'eor\`eme 1]{M}). As we can see in Figure~\ref{fgtact},
this interchange has the effect of taking the generators of 
$\pB = \dn = <a,b|a^2, b^3, (ab)^n>$ to their inverses.

Thus, if we take $\phi_0$ as the representation of $\dn$ induced 
by $\brhoz$, we have the following commutative diagram:

$$ \begin{CD}
\pt @>>> \pi_1(\St) @>>> \dn @>{\phi_0}>> \PSLC \\
@VV{{[} \alpha {]} \mapsto [\alpha]^{\mu}}V @V{\cong}VV 
@V{\tau_{\sharp}}VV \\
\pt @>>> \pi_1(\St) @>>> \dn @>{\phi_0}>> \PSLC 
\end{CD} $$

Since $\tau_{\sharp}$ takes $a$ and $b$ of $\dn$ to their inverses,
and $ab$ to $a^{-1}b^{-1}$ which is conjugate to $(ab)^{-1}$,
we see that 
$\mbox{tr}(\phi_0 \circ \tau_{\sharp}) = \mbox{tr}(\phi_0)$. On the other
hand, since $\brhoz$ is irreducible, $\phi_0$ is as well and 
we deduce (see \cite[Proposition 1.5.2]{CS1})
that there is an $A \in \PSLC$ with  
$\phi_0 \circ \tau_{\sharp} = A \phi_0 A^{-1}$. In other words, we have 
found an $A$ such that 
$A \brhoz(\beta) A^{-1} = \brhoz(\mu \beta \mu^{-1})$ for all $\beta
\in \pt$.

As mentioned above, we can complete the proof by showing that
$A^2 = \pm \I$.
Since $\tau$ is an involution, $\tau_{\sharp}^2 \equiv 1$ and 
$A^2$ commutes with every element of the irreducible 
representation $\phi_0(\dn)$. 
It follows that 
$A^2 = \pm \I$.  \qed

\section{Proof of Proposition~\ref{prpS}}

In this section we prove

\setcounter{section}{1}
\setcounter{theorem}{2}

\begin{prop} 
The minimum of the total norm $\| \cdot \|_T$ is $S = 3(|n-2| -1)$. 
\end{prop}

\setcounter{section}{3}
\setcounter{theorem}{0}

We break the proof down as a series of lemmas. 
Let us begin with an overview of the
strategy of the proof. 

Since $\mu$ is a cyclic surgery slope which is not a boundary
slope, $\| \mu \|_T = S$ (\cite[Corollary 1.1.4]{CGLS}) whence 

$$ S = \| \mu \|_T = 2 \| \mu \|_T - \| \mu \|_T = 
\| \mu^2 \|_T - \| \mu \|_T . $$

As the seminorms are given by the degree of $f_{\gamma}$, we
can determine $S$ by investigating the characters $x$ where the
order of 
zero differs: $Z_x(f_{\mu^2}) > Z_x(f_{\mu})$.
(By \cite[Proposition 1.1.3]{CGLS}, for all $x \in X$, $Z_x(f_{\mu^2}) \geq Z_x(f_{\mu})$.)  
That is,

$$S = \| \mu^2 \|_T - \| \mu \|_T = \sum_i \sum_{x \in \tXi} 
Z_x (f_{\mu^2}) - Z_x(f_{\mu}). $$

We argue that such a
``jumping point'' $x$
is the character of an irreducible representation $\rho$
(Lemma~\ref{lem41}).
Then $\rho$ is either a binary 
dihedral representation, or else $\brhoz$, 
the restriction to $\pt$ of the corresponding
$\PSLC$-representation, is non-abelian.
The number of irreducible dihedral characters 
(We will often refer to the characters of irreducible (dihedral, etc.) 
representations as irreducible (dihedral, etc.) characters.)
is easily related to 
the Alexander polynomial and we find that there are $(|n-6|-1)/2$ of
these (Lemma~\ref{lem42}). 

The non-abelian $\brhoz$'s will become representations of $\dn$.  
Ben Abdelghani and Boyer~\cite{BB} have 
calculated the number of characters of such a triangle group.
In the case of $\dn$, there are $(|n|-1)/2$ $\PSLC$-characters.
As each is covered twice, there are $|n|-1$ $\SLC$-characters 
(Lemma~\ref{lem43}).

In total then, we have $(|n-6| -1)/2  + |n|-1 = 3(|n-2|-1)/2$
(recall that $n$ is odd and not $1$, $3$, or $5$)
$\SLC$-characters $x$ where $Z_x(f_{\mu^2}) > Z_x(f_{\mu})$.
We finish the argument by showing that the difference in 
degree of zero is two at each of these points(Lemma~\ref{lem44}).

Let us now tackle the details of the argument.
Let $x \in \tXi$
with $Z_x(f_{\mu^2}) > Z_x(f_{\mu})$ where $X_i$ is an
algebraic component of $X$ containing an irreducible character.
Let $\nu : X_i^{\nu} \to X_i$ 
denote normalization (\cite[Chapter II, \S 5]{Sh}).
The birational map from  $\tXi$ to $X_i$ is regular at
all but a finite number of points of $\tXi$, called ideal points.
As in \cite[1.5]{CGLS}, $X_i^{\nu}$ may be identified with the complement
of the ideal points in $\tXi$.

\begin{lemma} \label{lem41} If $x \in \tXi$ with 
$Z_x(f_{\mu^2}) > Z_x(f_{\mu})$, then $\nu(x) = \chi_{\rho}$
is the character of an irreducible representation $\rho$.
\end{lemma}

\Pf
Suppose $x$ were an ideal point.
Since $Z_x(f_{\mu^2}) > Z_x(f_{\mu})$, 
$f_{\mu^2}(x) = 0$ so that $x$ is not a pole of $f_{\mu^2}$.
A little algebra shows 
\begin{equation} \label{eqT}
f_{\mu^2} = (f_{\mu})^2 + 4 f_{\mu}. 
\end{equation}
So, $x$ is also not a pole of $f_{\mu}$ 
and therefore $I_{\mu}(x) \neq \infty$ as well.
It follows that either $M$ admits a closed essential surface, or
else $\mu$ is a boundary slope (see \cite[Proposition 1.3.9]{CGLS}). 
However, since $K_n$ is a
Montesinos knot with less than four tangles, neither is true
(see~\cite[Section 1 and Corollary 4]{O}).

Thus, we can assume $x \in \Xin$ and 
write $\nu(x) = \chi_{\rho}$ with $\rho \in R_i$.
We wish to show that $\rho$ is an irreducible representation.
The idea is to show that, if $\rho$ is reducible, then the 
tangent space at $\rho$ is too small.

We can identify the Zariski tangent space
at $\rho$ with a subspace of the space of $1$-cocycles
$Z^1(\pi;\slC_{Ad\rho})$ (see \cite[Section 1.2]{G} or 
\cite[Section 3]{W}).
We can see that $R_i$ is four dimensional since, as we have mentioned,
$X_i$ is one dimensional (the knot being small) 
and by \cite[Corollary 1.5.3]{CS1}, 
$\dmn{R_i} = \dmn{X_i} + 3$.
Thus, $\dmn{Z^1(\pi;\slC_{Ad\rho})} \geq 4$.

Now, suppose $\rho$ were reducible. Then, by conjugating, 
we can take $\rho$ to be a representation into the upper triangular matrices.
Replace each matrix 
$\rho(g) = \mat{a}{b}{0}{a^{-1}}$
by 
$\rhod(g) = \mat{a}{0}{0}{a^{-1}}$
to obtain a diagonal representation $\rhod$ with the same character $\nu(x)$.

Since 
$$ \mat{1/n}{0}{0}{n} \mat{a}{b}{0}{a^{-1}} \mat{1/n}{0}{0}{n}^{-1} =  
\mat{a}{b/n^2}{0}{a^{-1}},$$
and $R_i$ is closed under conjugation \cite[Proposition 1.1.1]{CS1}, 
we can find representations on $R_i$ arbitrarily close to $\rhod$.
But, as $R_i$ is closed, $\rhod \in R_i$. So without loss of
generality, we can assume $\rho$ is diagonal.

\begin{claim} $\rho(\mu) = \mat{\pm i}{0}{0}{\mp i}$ and 
 $\rho(\pi) \cong \Z/4$. 
\end{claim}

\Pf (of Claim)
By, Equation~\ref{eqT},
$Z_x(f_{\mu^2}) > Z_x(f_{\mu})$ implies $Z_x(f_{\mu}) = 0$ 
and therefore, $\mbox{trace}(\rho(\mu)) \neq \pm 2$. On the
other hand, $\mbox{trace}(\rho(\mu^2))$ {\em is} $\pm 2$,
so $\rho(\mu^2) = \pm I$ (we're assuming that $\rho$ is diagonal).
Then $\rho(\mu) = \mat{\pm i}{0}{0}{\mp i}$. Since $\pi$ is 
normally generated by $\mu$, $\rho(\pi) \cong \Z/4$.
\qed (Claim)

Given this, we can calculate $\dmn{Z^1(\pi;\slC_{Ad\rho})}$ 
directly. Using~\cite[Theorem 1.1(i)]{BN}, 
$\dmn{H^1(\pi;\slC_{Ad\rho})} = b_1(\pi; \slC_{Ad\rho}) = 1$.
(This argument is explained in more detail and in a more general 
context in the proof of Lemma~\ref{lem51}.)
We can also determine $\dmn{B^1(\pi;\slC_{Ad\rho})}$ as
we have the surjection 
$$ \begin{array}{ccc} 
\slC & \longrightarrow & B^1(\pi;\slC_{Ad\rho}) \\
  A  & \mapsto & (u_A : \gamma \mapsto A - Ad\rho (\gamma)(A)). 
\end{array} $$ 
Since $\rho(\mu) = \mat{\pm i}{0}{0}{\mp i}$, the kernel 
is the one-dimensional set $\{ A \in \slC \mid A = \mat{a}{0}{0}{-a} \}$
while $\slC$ has dimension 3. Therefore, $\dmn{B^1(\pi;\slC_{Ad\rho})} =2$
and
\begin{eqnarray*}
\dmn{Z^1(\pi;\slC_{Ad\rho})} & = & \dmn{H^1(\pi;\slC_{Ad\rho})} + 
  \dmn{B^1(\pi;\slC_{Ad\rho})} \\
 & = & 1 + 2 = 3.
\end{eqnarray*}

This contradiction with our earlier estimate of the dimension of the 
cocycles shows that $\rho$ is irreducible. \qed

Thus, we can assume that $x$ is the character 
of an irreducible representation $\rho$. 
Then, by~\cite[Proposition 1.5.2]{CGLS}, $\rho(\mu^2) = \pm \I$
and $\brhoz$, 
the induced $\PSLC$-representation of $\pt$, will factor
through $\pd$. There are now two cases depending on whether
or not $\brhoz$ is abelian.

\begin{lemma} \label{lem42} 
If $\brhoz$ is abelian, then $\rho$ has (binary) dihedral 
image in $\SLC$. There are $(|n-6|-1)/2$ jumping points $x$ of this
type.
\end{lemma}

\Pf
If $\brhoz$ is abelian, it factors through
the finite group $H_1(\St; \Z)$. Thus, $\brhoz(\pt)$ is cyclic and
extending to $\pi$ and lifting we see that $\rho$ has binary
dihedral image in $\SLC$. 

Furthermore, any such dihedral representation {\em will} result in a jumping
point. For let $\nu(x) = \chi_{\rho}$ be the character of the 
binary dihedral $\SLC$-representation $\rho$. The corresponding
$\PSLC$ representation $\brho$ has as image a dihedral group normally
generated by $\brho(\mu)$. Therefore, $\brho(\mu)$ is of order two
and consequently $\rho(\mu) \neq \pm \I$ while $\rho(\mu^2) = \pm \I$.
This implies $\tr{\rho(\mu)} = 0$ so that $Z_x(f_{\mu^2}) > Z_x(f_{\mu})$.

The number of dihedral characters $d$ can be
related to the Alexander polynomial $\Ap$. Indeed, 
$d$ is equal to $(|\Apone|-1)/2$. (See \cite[Theorem 10]{K}
and recall that, up to conjugation, a dihedral subgroup of 
$\SLC$ may be assumed to lie in $SU(2)$.)

Finally, by \cite[Proposition 14]{Mr} or 
\cite[Theorem 1.2]{Hi}, 
$$\Ap \doteq (t-1)(t^{(n+3)}-1)/(t+1) + 
   t(t^3+1)(t^n+1)/(t+1)^2, $$
and the number of jumping points is $d = (|n-6|-1)/2$. \qed 

\begin{lemma} \label{lem43} 
If $\brhoz$ is non-abelian, then $\brhoz$ factors through to 
give an irreducible representation of $\dn$. There are
$|n|-1$ jumping points $x$ of this type.
\end{lemma}

\Pf
Since $\brhoz$ is a non-abelian $\PSLC$-representation
of $\pd$, it factors through to give an irreducible
representation of $\dn$, the fundamental group of the 
base orbifold of $\St$. (Recall that $\pd$ is a central 
extension of $\dn$.)

By \cite[Proposition D]{BB}, 
the number of $\PSLC$-characters of $\dpqr$ is
\begin{eqnarray} \label{eqpqr}
&& (p- \fltwo{p} - 1)(q- \fltwo{q} - 1)(r- \fltwo{r} - 1) +
   \fltwo{p} \fltwo{q} \fltwo{r} \\
&& \mbox{}+ \fltwo{\gcd(p,q)} + \fltwo{\gcd(p,r)} + \fltwo{\gcd(q,r)}
   + 1 \nonumber \end{eqnarray}
where $\fl{x}$ denotes the largest integer less than or equal to $x$.
This count includes the reducible characters. As the character of a 
reducible representation is also the character of an abelian representation,
we see that the reducible characters correspond to representations of
$H_1(\dpqr) = \Z / a \oplus \Z / (b/a)$ where $a = \gcd(p,q,r)$
and $b = \gcd(pq,pr,qr)$. So the number of reducible 
$\PSLC$-characters of $\dpqr$ is 
\begin{equation} \label{eqH1}
   \begin{array}{cl} \fltwo{b} + 1, & \mbox{if } a \equiv 1 \pmod{2}  \\
                     \fltwo{b} + 2, & \mbox{if } a \equiv 0 \pmod{2}.
   \end{array}
\end{equation}

Thus, there are $(|n| - 1)/2$ irreducible $\PSLC$-characters
of $\dn$. The corresponding representations each extend to an 
irreducible representation $\brhoz$ of $\pt$. These in turn 
can be extended to $\pi$ by Proposition~\ref{prpext}. Moreover,
as we see from the proof of that proposition, 
any representation $\brho$ which 
extends $\brhoz$ is such that $\brho(\mu)$ has order two. Thus, 
the irreducible
representations of $\dn$ all lead to jumping points where 
$Z_x(f_{\mu^2}) > Z_x(f_{\mu})$.
As these $\PSLC$-characters are covered twice in $\SLC$
(\cite[Lemma 5.5]{BZ}), we have $|n|-1$ $\SLC$ jumping points
of this type. \qed.

Combining the two previous lemmas, we have $3(|n-2| - 1)/2$
characters $\nu(x)$ where the degree of zero jumps.
We can complete the argument by showing that the jump is
two and that each $\nu(x)$ is a simple point of $X$ (i.e., 
$\nu(x)$ is on a unique irreducible component of $X$ and smooth on 
that component, see \cite[Chapter 2 \S 2]{Sh}).

\begin{lemma} \label{lem44}
At each jumping point $x$, $\nu(x)$ is a simple point
of $X$, and the jump is two, 
$$Z_x(f_{\mu^2}) - Z_x(f_{\mu}) = 2.$$
\end{lemma}

\Pf
We have argued that $\nu(x) = \chi_{\rho}$ is the 
character of an irreducible representation $\rho$. Moreover,
either $\rho$ is dihedral, or else $\brhoz$ factors through 
$\dn$. Dihedral representations of $2$-bridge knots 
have been analyzed by Tanguay~\cite{Ta} and we will adopt
his arguments to our situation. 

For representations
going through $\dn$, we plan to follow the argument 
of~\cite[Section 4]{BZ}
(see also \cite[Theorem A]{BB}). 
The essential requirements in this approach are that
$\rho(\pdM) \not\subset \{ \pm \I \}$, and that
$\nu(x)$ is a simple point of $X$.
However, as Claim~\ref{clmoct} below suggests, we will be
frustrated if $\rho$ is octahedral. 

Therefore, we will break the problem into 
three cases: $\brhoz$ factors through the triangle group
$\dn$ and 
$\rho$ is not octahedral; $\rho$ is dihedral; and $\rho$ is octahedral.

\bigskip

\mbox{ }

\bigskip

\noindent{\bf Case 1:} 
$\brhoz$ factors through $\dn$ and 
$\rho$ is not octahedral

\medskip

To apply the method of \cite[Section 4]{BZ}, we need to 
show that $\rho(\pdM) \not\subset \{ \pm \I \}$. The following
claim starts us on the way.

\begin{claim} \label{clmoct}
Suppose $\brhoz$ factors through $\dn$. 
If $\brhoz(\dpt) = \{ \pm \I \}$
then $\rho$ is a (binary) octahedral representation.
\end{claim}

\Pf (of Claim)
We've argued that $\brhoz(\mt) = \brhoz(\mu^2) = \pm I$, so we
are lead to investigate the image of $\lt$, the lift of $\lambda$
to $\pt$. Trotter~\cite{Tr} has explained how to find 
the image of $\lt$ in $\dpqr$ in the case of a $(p,q,r)$ pretzel
knot with $p,q,r$ all odd. Following the analogous procedure in 
our case, we find that $\lt$ projects to $ c^k a c^{k+1} b c b$ in  
$\dn = <a,b,c|a^2, b^3, c^{|n|}, abc^{-1}>$, where $k = \fltwo{ |n|}$.

We can take 
$$ \brhoz(c) = \pm  \left( \begin{array}{cc} \omega & 0  \\
                                0 & {\omega}^{-1} \end{array} \right)
\mbox{ , } 
   \brhoz(b) = \pm  \left( \begin{array}{cc} u & 1  \\
                           u(1-u) - 1  & 1-u \end{array} \right)$$
where $\omega = e^{\pi i j / |n|}$, $ 1 \leq j \leq \fltwo{|n|}$
and $u \in \C$ (compare \cite[Example 3.2]{BZ2}).
Since $a = cb^{-1}$ is of order two, $\tr{\brhoz(a)} = 0$
whence $u = \omega^2/(\omega^2 - 1)$.

Now, 
$$\brhoz(c^k a c^{k+1}) = \brhoz(c^{k+1}b^{-1}c^{k+1})
= \Pmat{(-1)^j (1-u) \omega}{-1}{1+u(u-1)}{(-1)^ju/ \omega}$$
while 
$$\brhoz(bcb) = \Pmat{\omega u^2 + (u-u^2-1)/{\omega}}
{\omega u + (1-u)/{\omega}} {(u^2-u^3-u) \omega + \frac{u^3-2u^2+2u-1}{\omega}}
{(u-u^2-1) \omega + \frac{(1-u)^2}{\omega}}$$
and so, after substituting $u = \omega^2/(\omega^2-1)$,
\begin{eqnarray*}
\tr{\brhoz}(\lt) &=& \pm 2(\frac{\omega}{\omega^2-1})^2
{[} \frac{\omega^6+1}{\omega^3} - (-1)^j \frac{\omega^4+1}{\omega^2} {]} \\
& = & \pm  \frac{1}{\sin^2 (\pi j/|n|)}{[} \cos (3 \pi j /|n|) - (-1)^j
\cos(2 \pi j / |n|) {]} 
\end{eqnarray*}
Thus, $\tr{\brhoz}(\lt) = \pm 2$ only if $\pm \sin^2(\pi j/|n|)$ is 
$\cos(5 \pi j/2 |n|) \cos(\pi j/ 2|n|)$ or
$\sin(5 \pi j/2 |n|) \sin(\pi j/ 2|n|)$. The only choice consistent with
our conditions on $j$ and $n$ is that $j/|n| = 1/3$. 

In other words,
as long as $j/|n| \neq 1/3$, we are assured that
$\tr{\brhoz}(\lt) \neq \pm 2$ and, therefore, that 
$\brhoz(\dpt) \neq \{ \pm I \}$. 

On the other hand, if $j/|n| = 1/3$ 
(this is possible only when $3 \mid n$), then
$\rho$ is a (binary) octahedral representation. \qed (Claim)

By Claim~\ref{clmoct}, $\brhoz(\dpt) \neq \{ \pm I \}$.
Then, $\rhoz(\dpt) \not\subset \{ \pm \I \}$ and
$\rho(\pdM) \not\subset \{ \pm \I \}$ as well.

The other requirement of \cite[Section 4]{BZ}
is that $\nu(x)$ be simple. 
We first show that the corresponding character 
$y =\chi_{\rhoz}$ is smooth in $Y$, the character variety for $\pt$.
As the Zariski tangent space at $\rhoz$ can be identified
with a subspace of the cocycles, we proceed by investigating the group
cohomology.

\begin{claim}
$Z^1(\pd; {\slC}_{Ad \brhoz}) \cong
Z^1(\dn; {\slC}_{Ad \brhoz})$. 
\end{claim}

\Pf (of Claim)
The Seifert structure of $\St$
gives the exact sequence
$$ 0 \longrightarrow F \longrightarrow \pd \stackrel{\phi}{\longrightarrow}
\dn \longrightarrow 1$$
where $F = \langle h \rangle \cong \Z$ is the group of a regular fibre $h$.
The projection $\phi$ induces a homomorphism $\Phi: 
Z^1(\dn; {\slC}_{Ad \brhoz}) \to Z^1(\pd; {\slC}_{Ad \brhoz})$.

To construct the inverse, we show that, for
each $u \in Z^1(\pd;{\slC}_{Ad \brhoz})$, $u(h) = 0$. Indeed,
for all $g \in \pd$, $u(hg) = u(gh)$. On the other hand, 
$\brhoz(h)$ commutes with $\brhoz(\dn)$. Since $\brhoz$ is irreducible,
this implies $\brhoz(h) = \pm \I$. Putting it together,
\begin{eqnarray*}
u(h) + u(g) & = & u(h) + Ad \brhoz(h) \cdot u(g) \\
            & = & u(hg) \\
            & = & u(gh) \\
            & = & u(g) + Ad \brhoz(g) \cdot u(h) \\
            & = & u(g) + \brhoz(g) u(h) \brhoz(g)^{-1}.
\end{eqnarray*}

Thus, $u(h) \in \slC$ commutes with $\brhoz(\pd)$ which implies $u(h) = 0$.
Now define $\Psi :
Z^1(\pd; {\slC}_{Ad \rhoz}) \to Z^1(\dn; {\slC}_{Ad \rhoz})$
by $\Psi(u)(\phi(g)) = u(g)$. Since $u(h) = 0$, $\Psi$ is well defined.
Moreover, it's an inverse of $\Phi$ and we have the required isomorphism.
\qed (Claim)

This cocycle isomorphism also descends to the level of cohomology:
$$H^1(\pd; {\slC}_{Ad \brhoz}) \cong H^1(\dn; {\slC}_{Ad \brhoz}).$$
On the other hand, a straightforward calculation shows that 
the cohomology of the triangle group is trivial (for example,
see \cite[Lemma~5.1.3]{Ma1}).

Now, since the $\PSLC$-representation $\brhoz$ and the 
$\SLC$-representation $\rhoz$ result in exactly the same
adjoint action on $\slC$, we see that 
$$\dimC{H^1(\pd; {\slC}_{Ad \rhoz})} = 0$$ as well.
So we can proceed as in~\cite[Section 4]{BZ} to show that 
$H^1(\pt;{\slC}_{Ad \rhoz})$ has dimension one
and hence that $y$ is simple in $Y$.

The following lemma allows us to relate the smoothness of 
$y$ to that of $\nu(x)$.

\begin{lemma}
Let $\rho$ be an $\SLC$-representation of
a finitely generated group $\pi$ and $\rhoz$ the restriction
to a normal subgroup of finite index $\pt$. 
Then 
$$\mbox{dim}_{\C}H^1(\pi;{\slC}_{Ad \rho}) 
\leq \mbox{dim}_{\C}H^1(\pt;{\slC}_{Ad \rhoz}).$$
\end{lemma}

\Pf
The Lyndon - Hochschild - Serre spectral sequence gives us the
exact sequence (see~\cite[Theorem 11.5]{Rt})
$$ 0 \longrightarrow H^1(\pi/\pt; ({\slC}_{Ad \rho})^{\pt}) \longrightarrow 
   H^1(\pi; \slC_{Ad \rho}) 
\longrightarrow H^1(\pt; \slC_{Ad \rhoz})^{\pi/\pt},$$
where $A^G = \{a \in A \mid g \cdot a = g, \forall g \in G \}$ 
denotes the set of fixed points of the module $A$ under the
group action $G$.
Now, $H^1(\pi/\pt; (\slC_{Ad \rho})^{\pt}) = 0$ since
$\pi/\pt$ is finite and $(\slC_{Ad \rho})^{\pt}$ is a complex
vector space.
So we have
$$ H^1(\pi; \slC_{Ad \rho}) \hookrightarrow 
  H^1(\pt; \slC_{Ad \rhoz})^{\pi/\pt}
  \hookrightarrow H^1(\pt; \slC_{Ad \rhoz}).$$ \qed

In our case, the proposition shows that 
$\mbox{dim}_{\C}H^1(\pi;\slC_{Ad \rho}) \leq 1$ whence $\nu(x)$ is a 
smooth point of $X_i$ and a simple point of $X$.

\medskip

\Rmk We have been using the ideas of~\cite[Section 4]{BZ} whereby,
under appropriate conditions, $x = \chi_{\rho}$ is simple in $X(\pi)$
exactly when $\mbox{dim}_{\C}H^1(\pi;\slC_{Ad \rho}) = 1.$
Interpreted in this context, the proposition says
``simple points of $X(\pt)$ lift to simple points of $X(\pi).$''

\medskip

Thus, if $\rho$ factors through $\dn$ and is not octahedral,
then $\nu(x)$ is a simple point of $X$ 
and $\rho(\pdM) \not\subset \{ \pm \I \}$. Following the reasoning
of~\cite[Section 4]{BZ} 
(see also \cite[Theorem A]{BB}), 
we conclude that $Z_x(f_{\mu^2}) - Z_x(f_{\mu}) = 2$.
This proves Lemma~\ref{lem44} in the case
that 
$\brhoz$ factors through $\dn$ and 
$\rho$ is not octahedral.

\medskip

\noindent{\bf Case 2:} $\rho$ is dihedral

\medskip

Here we follow the approach outlined by Tanguay in his 
thesis~\cite{Ta}. Since that document remains unpublished, we
give a summary of his argument. As in the previous case, we
appeal to \cite[Section 4]{BZ} (or 
\cite[Theorem A]{BB}). Most of the work goes into proving
that $H^1(\pi;\slC_{Ad \rho})$ has dimension $1$ when $\rho$ is 
a dihedral representation (Claim~\ref{clmdi}). Once this is established, 
\cite[Lemma 4.5]{BZ} shows that $\nu(x)$ is a simple point of 
$X$. Moreover, as the non-abelian group $\rho(\pi)$
is generated by $\rho(\pdM)$, it follows that 
$\rho(\pdM) \not\subset \{ \pm I \}$. Lemmas 4.6 through 4.9 of
\cite{BZ} then imply $Z_x(f_{\gamma}) = 2$ for each non-trivial 
$\gamma \in L$ such that $f_{\gamma}(x) = 0$. In particular, we've already
argued that $\rho(\mu^2) = \pm I$ (see proof of Lemma~\ref{lem42}) whence
$f_{\mu^2}(x) = 0$, and, therefore, $Z_x(f_{\mu^2}) = 2$. 
On the other hand,
$\tr{\rho(\mu)} = 0$ (proof of Lemma~\ref{lem42}) implies
$Z_x(f_{\mu}) = 0$. Thus, $\nu(x)$ is simple and the jump at $x$,
$Z_x(f_{\mu^2}) - Z_x(f_{\mu})$, is two.

Therefore, the following claim will suffice to prove Lemma~\ref{lem44}
in the case that $\rho$ is dihedral.

\begin{claim} \label{clmdi}
$\dim{H^1(\pi;\slC_{Ad \rho})} = 1$ for $\rho$ dihedral.
\end{claim} 

\Pf (of Claim)
In \cite[Theorem 1.1(i)]{BN}, Boyer and Nicas 
explain how the dimension of the cohomology with coefficients twisted
by a cyclic representation can be written in terms of 
certain covering spaces of the manifold.
The plan is to adapt that argument to the present
situation of a dihedral representation.

Let $\rho(\pi) = D_{4m}$, the binary 
dihedral group of order $4m$. 
Then $Ad \rho (\pi) \subset \mbox{Aut}(\SLC)$ is isomorphic to
$D_{2m}$, the dihedral group of order $2m$.
In analogy with \cite[Theorem 1.1(i)]{BN}, 
Tanguay~\cite{Ta} shows that the Betti number
$b_1(\pi; \slC_{Ad \rho})$ can be related to the Betti 
numbers of several covers of $M$:
$$ b_1(\pi; \slC_{Ad \rho}) = b_1 (\pt ; \C) - b_1 (\pi; \C)
 + \frac{1}{\phi(m)} \sum_{d|m} \mu(\frac{m}{d}) b_1( \pi_d; \C),$$
where the $\pi_d$ are the kernels of the maps
$$\pi \stackrel{Ad \rho}{\longrightarrow} D_{2m} \to  D_{2d},$$
and $\phi$ and $\mu$ are the Euler and M\"obius Functions
respectively. Now, $b_1(\pi; \C)=1$~\cite[Exercise 2.E.6]{R} and
$b_1(\pt; \C) = 1$~\cite[Section 8.D]{R}. So we can show that 
$b_1(\pi; \slC_{Ad \rho}) =1$ 
by arguing that $b_1(\pi_d; \C) = d$.

Let $\Md$ be the covering of $M$ corresponding to $\pi_d$.
Then $\Md$ also covers $\Mt$ and this covering may be extended
to an orbifold covering $\Sd \to \St$. We will argue that 
$b_1(\Sd) = 0$. Then, since $\Sd$ is obtained from $\Md$
by filling along $d$ tori, $0 = b_1(\Sd) \geq b_1(\Md) - d$
whence $b_1(\Md) \leq d$. On the other hand, as 
$\Md$ has $d$ toral boundary components, Lefschetz duality
allows us to argue that $b_1(\pi_d) = \dmn{H_1(\Md ; \C)} \geq d$.
Therefore, $b_1(\pi_d) = d$, as required.

It remains to show that $b_1(\Sd) = 0$, and this is where we
must introduce some new ideas beyond those used by Tanguay. 
We have the diagram
$$
\begin{CD}
0 @>>> E @>>> \pi_1(\Sd) @>>> \pi_1^{\mbox{orb}}(B_d) @>>> 1 \\
@VVV @VV P' V @VV P V @VV P''V \\
0 @>>> F @>>> \pd @>>> \dn @>>> 1
\end{CD}
$$
where the horizontal rows are the exact sequences arising 
from the Seifert structure of $\Sd$ and $\St$, $E \cong F \cong \Z$
represent regular fibres, and $B_d$ is the base orbifold
of $\Sd$. Now, $\img{P}$ is normal since $P$ is a regular
cover. This implies $\img{P''}$ is normal. Since $E$ and $F$
are abelian, $\img{P'}$ is also normal. Thus, the cokernels will
be groups and we can use the 
Snake Lemma to obtain the exact sequence
$$ \noy{P''} \stackrel{\delta}{\longrightarrow} \coker{P'}
\stackrel{\alpha}{\longrightarrow} \coker{P}
\stackrel{\beta}{\longrightarrow} \coker{P''} \longrightarrow 1. $$
Here, $\noy{P''} = 0$ since $B_d \to \dn$ is an orbifold covering space.
Thus, $\alpha$ is injective.
Since $P$ comes from the dihedral covering $\Md \to \Mt \to M$, we see 
that $\coker{P} \cong \Z / d$. By the injectivity of $\alpha$,
$\coker{P'} \cong \Z / a$ where $a \mid  d$. On the other
hand, as the degree $d$ of the Seifert cover $\Sd \to \St$ is 
the product of the degree $a$ in the fibres and the degree of
the orbifold cover $c$, we see that $\card{\coker{P''}} = c = d/a$.
However, since $\img{\al} = \noy{\beta}$, 
$\img{\beta}$ also has cardinality $c$:
$$\coker{P''} = \img{\beta} \cong \coker{P} / \noy{\beta}
\cong (\Z /d) / (\Z / a) \cong \Z / c.$$ 

The projection $\dn \to \coker{P''} \cong \Z / c$ 
is an abelian representation, and as such
factors through $H_1(\dn) = \Z / b$, where 
$b = \mbox{gcd}(3n,6,2n) = \mbox{gcd}(3,n) = 1 \mbox{ or } 3$. So the covering 
$B_d \to S^2(2,3,|n|)$ is either trivial, or else of degree $3$
whence $B_d$ is either $S^2(2,3,|n|)$ or
else $S^2(2,2,2,3, |n|/3)$.
Given $B_d$, we have an explicit formula (see \cite[Equation 12.31]{BuZ})
for $\pi_1(\Sd)$ involving the orders of the cone points. We can
then show that $H_1(\Sd; \Z)$ is torsion by examining its order ideal
(see~\cite[Section 8.B]{R}). Therefore, $b_1(\Sd) = 0$ as required.

This completes the proof of the claim. \qed (Claim)

\medskip

\noindent{\bf Case 3:} $\rho$ is octahedral

\medskip

The approach here is much like that used for the dihedral 
representations. In this case, 
$$b_1( \pi; \slC_{Ad \rho}) = b_1(\hp; \C) - b_1(\pi; \C)$$
where $\hp = \rho^{-1}(D_6)$
with $D_6$ a dihedral subgroup of index four in the octahedral
group $S_4$.
Of course $b_1(\pi; \C) = 1$ as before, so we will want to argue
that $b_1(\hp; \C) = 2$.

Let $\hM$ be the covering of $M$ corresponding to $D_6$. 
Then $\hM$ is an irregular covering of degree $4$ which
also covers $\Mt$. As before, we extend the covering
$\hM \to \Mt$ to a degree two map between Seifert spaces: $\hS \to \St$.
This leads to a diagram quite similar to that for the dihedral
representation:
$$
\begin{CD}
0 @>>> E @>>> \pi_1(\hS) @>>> \pi_1^{\mbox{orb}}(\hB) @>>> 1 \\
@VVV @VV P' V @VV P V @VV P''V \\
0 @>>> F @>>> \pd @>>> \dn @>>> 1
\end{CD}
$$

In this case, $\hB$, the base orbifold of $\hS$, is either a $1-1$
or a $2-1$ cover of $S^2(2,3,|n|)$.
In particular it is a regular cover and $\coker{P''}$ is either
trivial or cyclic of order $2$. Since $\coker{P''}$ is abelian,
it is also a factor of $H_1(\dn) \cong \Z /b$
where $b = \mbox{gcd}(3,n) = 1 \mbox{ or } 3$. Thus, $\hB = S^2(2,3,|n|)$
as well and, as in the dihedral case, we find $b_1(\hS) = 0$.
Since $\hM$ has two boundary components, we may now argue that
$b_1(\hp; \C) = 2$. Thus, $b_1( \pi; \slC_{Ad \rho}) = 1$ and
$\nu(x)$ is again a smooth point of $X_i$ (and
a simple point of $X$) yielding a jump of two.

This completes the proof of Lemma~\ref{lem44} and with it,
the proof of Propostion~\ref{prpS}. \qed

\section{Proof of Proposition~\ref{prp2n4}}

In this section, we prove

\setcounter{section}{1}
\setcounter{theorem}{3}

\begin{prop} 
The total norm of the $2n+4$ Seifert fibred surgery is
$\| 2n+4 \|_T = S + 3(|n-6| - 1)$.
\end{prop}

\setcounter{section}{4}
\setcounter{theorem}{0}

\Pf
Bleiler and Hodgson~\cite[Proposition 17]{BH} 
have shown that $2n+4$
surgery on the $(-2,3,n)$ pretzel knot
$K_n$ results in a manifold which 
is Seifert fibred over $S^2(2,4,|n-6|)$.
(Actually, there is a small error in the statement of 
their proposition which refers to ``$4n+14$ surgery on 
the $(-2,3,2n+7)$ pretzel knot.'' It should read ``4n+18
surgery on the $(-2,3,2n+7)$ pretzel knot.'')  
We can find the total Culler-Shalen norm
of this slope in much the same way as we calculated $S$ above.

Recall that
$$\| 2n + 4 \|_i = \sum_{x \in \tXi} Z_x(f_{ 2n + 4}) $$
where $Z_x(f)$ denotes the order of zero of $f$ at $x$.
Since the meridian $\mu$ of $K_n$ is not a boundary slope,
$Z_x(f_{\mu}) \leq  Z_x(f_{2n+4})$ for each $x$ 
(\cite[Proposition 1.1.3]{CGLS}).
This suggests that we approach the calculation of 
the total norm $\| 2n+4 \|_T$ by comparison
with $\| \mu \|_T = S:$
\begin{equation} \label{eqa1}
\| 2n+4 \|_T = S + \sum_i \sum_{x \in \tXi}( Z_x(f_{ 2n+4 }) - Z_x(f_{\mu})).
\end{equation}

We first show that the ``jumping points'' $x$,
where $Z_x(f_{\mu}) < Z_x(f_{2n+4})$, are characters of 
irreducible representations (Lemma~\ref{lem51}). We next 
show that there are $3(|n-6| -1)/2$ such characters
(Lemma~\ref{lem52}) and finally that the jump at each 
such character is two (\cite[Theorem A]{BB}). 

\begin{lemma} \label{lem51}
If $x \in \tXi$ with 
$Z_x(f_{2n+4}) > Z_x(f_{\mu})$, then $\nu(x) = \chi_{\rho}$
is the character of an irreducible representation $\rho$.
\end{lemma}

\Pf
Since $M$ is small and $2n+4$ is not
a strict boundary class, we may apply \cite[Proposition 1.6.1]{CGLS}
to see that $Z_x(f_{2n+4}) = Z_x(f_{\mu})$ at ideal points. So
we can assume $x \in X_i^{\nu}$. Then $\nu(x) = \chi_{\rho}$
is the character of a representation $\rho \in R_i$. 
We wish to show that $\rho$ is an irreducible representation.
The idea is to show that, if $\rho$ is reducible, then the 
tangent space (which is a subspace of the space of 
$1$-cocycles of $\pi$ by \cite[Section 1.2]{G} or 
\cite[Section 3]{W}) is too small.

As in the proof of Lemma~\ref{lem41}, 
we can argue that the dimension of $Z^1(\pi;\slC_{Ad\rho})$ 
is at least $4$. 

On the other hand, if $\rho$ were reducible, then, 
since $R_i$ is closed and invariant under conjugation, 
we can assume that $\rho$ is diagonal. Now,
$Z_x(f_{2n+4}) >  Z_x(f_{\mu})$
implies $\rho(2n+4) = \pm I$ (\cite[Proposition 1.5.4]{CGLS}).
So, if we take $\brho$ as the $\PSLC$-representation corresponding
to $\rho$, then $\brho$ factors through $H_1(M(2n+4); \Z) \cong \Z / (2n+4)$.
Since $\rho(\pi)$ is normally generated by $\rho(\mu)$ and
$\rho$ is diagonal, we see that $\rho(\mu) = \mat{\eta}{0}{0}{\eta^{-1}}$
with $\eta^{2(2n+4)} = 1$.

\begin{claim} $b_1(\pi; \slC_{Ad\rho}) = 1$.
\end{claim}

\Pf (of Claim)
By \cite[Theorem 1.1(i)]{BN},
$$b_1(\pi; \slC_{Ad\rho}) = b_1(\pi; \C) + 2 b_1(\pi; \C_{\beta})$$
where $\beta = \eta^2$ is a $(2n+4)$th root of unity. Of course 
\cite[Exercise 2.E.6]{R}, 
$b_1(\pi; \C) = 1$. Now, $\C_{\beta}$ is $\C$ with the $\Z$-action
induced by $t \cdot c = c \beta$ where $t$ is a generator of $\Z$.
Since $\pi$ surjects onto $H_1(M; \Z) \cong \Z$, this gives a $\pi$-action
on $\C$. 

We can also think of $H_1(M; \Z)$ as acting on $\hM$, the infinite cyclic
cover of $M$, and define a $\C{[}t, t^{-1}{]}$-module structure on 
$H_1(\hM; \C)$ (see \cite[Section 7.A]{R}). In this context
$H^1(\pi; \C_{\beta}) = \coker{H^1(\hM ; \C) 
\stackrel{t-\beta}{\longrightarrow} H^1(\hM; \C)}$ where
$t-\beta$ represents multiplication by $t - \beta$.
Since the Alexander polynomial $\Apn$ is the generator of $H_1(\hM; \C)$ 
as a $\C{[}t, t^{-1} {]}$-module, we can argue that
$\coker{t- \beta} = 0$ if $\Delta_{K_n}(\beta) \neq 0$.

It is not difficult to show that $\Apn$ admits no roots
which are $2n+4$th roots of unity (see Lemma~\ref{leArts}). 
Thus $H^1(\pi; \C_{\beta}) = 0$
and $b_1(\pi; \slC_{Ad\rho}) = 1$. \qed (Claim)

We can now argue as in 
the proof of Lemma~\ref{lem41} that $\dmn{B^1(\pi; \slC_{Ad \rho})} = 2$
and $\dmn{Z^1(\pi; \slC_{Ad \rho})} = 3$. This contradicts
our earlier estimate for the dimension of the cocycles
and we conclude that $\rho$ is an irreducible representation. 
\qed 

\begin{lemma} \label{lem52} 
There are $3(|n-6|-1)/2$ jumping points $x$ 
where $Z_x(f_{2n+4}) > Z_x(f_{\mu})$.
Moreover, at each jumping point, $\nu(x)$ is a simple point of $X$. 
\end{lemma}

\Pf
The plan is to argue that such a jumping point $\nu(x) = \chi_{\rho}$ 
gives rise to a $\PSLC$-representation $\brho$ which factors 
through $\Delta(2,4,|n-6|)$. Conversely, every such $\brho$
factoring through $\Delta(2,4,|n-6|)$ 
leads to a jumping point $\nu(x)$. Since the $\nu(x)$ 
are simple (and therefore each lies on only one curve in $X$), 
we can count the number of jumping points simply by counting
the number of $\PSLC$-representations of $\Delta(2,4,|n-6|)$.
Here are the details:

By \cite[Proposition 1.5.4]{CGLS}, $Z_x(f_{2n+4}) >  Z_x(f_{\mu})$
implies $\rho(2n+4) = \pm I$. Therefore, the corresponding 
$\PSLC$ representation $\brho$ factors through $\pi_1(M(2n+4))$.
As this is an irreducible $\PSLC$-representation (and therefore either
abelian or else with image $\Z/2 \times \Z/2$) it must
kill the center of $\pi_1(M(2n+4))$ and factor through
to give an irreducible $\PSLC$-representation $\brho'$ of 
$\Delta(2,4,|n-6|)$, the base orbifold of $M(2n+4)$.

Now, as in the proof of Lemma~\ref{lem44},
$$H^1(\pi_1(M(2n+4));\slC_{Ad \rho})
\cong H^1(\Delta(2,4,|n-6|);\slC_{Ad \brho'})$$ 
is trivial. Thus, arguing as in
\cite[Section 4]{BZ}, we can deduce that $\nu(x)$ is a smooth point of 
$X_i$ (and in fact a simple point of $X$)
so that $\nu^{-1}(\nu(x)) = x$.
Therefore, the jumping points $\nu(x) = \chi_{\rho}$
where $Z_x(f_{2n+4}) > Z_x(f_{\mu})$
are simple points of $X$ and correspond to irreducible
$\PSLC$ characters $\brho$
which factor through $\Delta(2,4,|n-6|)$. 

Conversely, any such 
representation induces a jumping point. This is immediate if 
the representation is diagonalizable on $\pdM$ since then
$\rho(\mu)$ is of finite order, but not $\pm \I$. On the other hand,
$\rho(2n+4) = \pm \I$. Thus $Z_x(f_{2n+4}) > 0 = Z_x(f_{\mu})$.
If $\rho(\pdM)$ is parabolic, we can appeal to \cite[Theorem A]{BB}. 

So, to find the number of jumping points, we must count the irreducible 
$\PSLC$-characters of $\Delta(2,4,|n-6|)$. 
By Equations~\ref{eqpqr} and \ref{eqH1},
there are $|n-6| - 1$ such. 
Since $d = (\Apone - 1)/2 = (|n-6| -1)/2$ (see Proof of Lemma~\ref{lem42}), 
half of these are dihedral characters.
Now, dihedral characters are covered once in $\SLC$ while
other characters are covered twice 
(\cite[Lemma 5.5]{BZ}). Thus we have $3 (|n-6|-1)/2$ $\SLC$-characters 
coming from irreducible $\PSLC$-representations 
which factor through $\Delta(2,4,|n-6|)$.
As we have argued, this is the number of jumping points where
$Z_x(f_{2n+4}) > Z_x(f_{\mu})$. \qed

To complete the proof of Proposition~\ref{prp2n4}, we use
\cite[Theorem A]{BB} which says that the jump is two at each 
jumping point. So in Equation~\ref{eqa1}, we have $3 (|n-6|-1)/2$
jumping points, each providing a jump of two. We conclude that
$\| 2n+4 \|_T = S + 3(|n-6| -1)$ and 
Proposition~\ref{prp2n4} is proved. \qed 

\section{Proof of Proposition~\ref{prp2n5}}

In this section, we prove 

\setcounter{section}{1}
\setcounter{theorem}{4}

\begin{prop} 
The total norm of the $2n+5$ Seifert fibred surgery is
$\| 2n+5 \|_T = S + 4(|n-5| - 2)$.
\end{prop}

\setcounter{section}{5}
\setcounter{theorem}{0}

\Pf
Bleiler and Hodgson~\cite[Proposition 16]{BH} 
have shown that $2n+5$
surgery on the $(-2,3,n)$ pretzel knot
$K_n$ results in a manifold which 
is Seifert fibred over $S^2(3,5,|n-5|/2)$ 
The argument here is essentially identical to that in
the previous section. In particular,
\begin{equation} \label{eqa2}
\| 2n+5 \|_T = S + \sum_i \sum_{x \in \tXi}( Z_x(f_{ 2n+5 }) - Z_x(f_{\mu})).
\end{equation}
The main difference is that now $x$ may be the character of
a reducible representation. The argument of Lemma~\ref{lem51}
depended on the Alexander polynomial having no zeroes at
a $(2n+4)$th root of unity. However, it may have a zero at a 
$(2n+5)$th root of unity. To clarify this point, we begin by 
investigating the roots of the Alexander polynomial.

\begin{lemma} \label{leArts}
Let $\Apn$ be the Alexander polynomial of the
$(-2,3,n)$ pretzel knot $K_n$. (In particular, $n$ is odd.)
 Suppose $\Delta_{K_n}(\zeta) = 0$
where $\zeta$ is a primitive $m$th root of unity. Then one
of the following is true.
\begin{itemize}
\item $3 \mid n$ and $m= 6$.
\item $10 \mid (n-1)$ and $m=10$.
\item $12 \mid (n-3)$ and $m=12$.
\item $15 \mid (n-5)$ and $m=15$.
\end{itemize}
\end{lemma}

\Pf
The Alexander polynomial never admits zeroes which are
prime power roots of unity. Indeed,
by~\cite[Theorem 8.21]{BuZ}, $H_1(\Sigma_m; \Z)$ is finite 
iff no root of the Alexander polynomial is an $m$th
root of unity. Here $\Sigma_m$ denotes the $m$-fold
branched cyclic cover of the knot (see \cite[Section 10.C]{R}). 
Using the Milnor~\cite{Mi}
sequence we can show that $b_1(\Sigma_m) = 0$ whenever
$m$ is a prime power. 

It is straightforward to argue that there are no zeroes at
$m$th roots of unity when $m \geq 18$. This leaves
$m = 6,10,12,14,$ and $15$. Since 
the value of $\Delta_{K_n}( \zeta )$ depends only on the value of
$n \bmod m$, one need only make the calculation for
each of the $m$ equivalence classes to verify that 
there are no roots when $m = 14$, and roots in the 
other cases ($m = 6, 10,12,15$) only as given in
the statement of the lemma. \qed

In particular, if $\Apn$ admits a $2n+5$ root of 
unity, then $m \mid 2n+5$ for one of $m =6,10,12 \mbox{ or } 15$.
However, if $6 \mid 2n+5$, then $3 \nmid n$.
Similarly, $m = 10$ and $m=12$ are not feasible. However,
if $m = 15$, we must have $15 \mid 2n+5$ and 
$15 \mid (n-5)$ which implies $n \equiv 5 \pmod{30}$. 
Indeed, one can verify that when $n \equiv 5 \pmod{30}$,
the Alexander polynomial admits primitive $15$th roots of unity as zeroes
and that these are simple (i.e., not repeated) zeroes.

So, we are lead to examine two cases.

\medskip

\noindent{\bf Case 1:} $n \not\equiv 5 \pmod{30}$

\medskip

In this case, the argument is identical to that of 
Proposition~\ref{prp2n4}, so we will omit most of the 
details. The jumping points correspond to irreducible $\PSLC$-characters
of $\Delta(3,5, |n-5|/2)$, the base orbifold of $M(2n+5)$.
By Equations~\ref{eqpqr} and \ref{eqH1}, there
are $|n-5|-2$ such. Since $H_1(M(2n+5); \Z) \cong \Z/(2n+5)$ has odd
order, none of these are dihedral characters. Therefore, they
are covered twice in $\SLC$ (\cite[Lemma 5.5]{BZ}). Since each contributes
a jump of two (\cite[Theorem A]{BB}), we have proved the proposition
in this case.

\medskip

\noindent{\bf Case 2:} $n \equiv 5 \pmod{30}$

\medskip

As in the proof of Lemma~\ref{lem51},  
we can apply \cite[Proposition 1.6.1]{CGLS}
to deduce that at a jumping point $x$ (i.e., where 
$Z_x(f(2n+5)) > Z_x(f_{\mu})$), $\nu(x) = \chi_{\rho}$
is the character of a representation $\rho$. The difficulty
is that now $\rho$ may be reducible.
We will consider irreducible and reducible jumping points separately.

\begin{lemma} \label{lem62}
 There are $2(|n-5|-6)$ jumping points where 
$\nu(x)$ is the character of an irreducible representation.
These are simple points of $X$ and each contributes a jump of two.
\end{lemma}

\Pf
If $\nu(x) = \chi_{\rho}$ is the character of an irreducible 
representation, we can follow the arguments of the proof of 
Lemma~\ref{lem52} to see that
$x$ is a simple point of $X$, and that the corresponding 
$\PSLC$-representation $\brho$ factors through $\Delta(3,5,|n-5|/2)$,
the group of the base orbifold of $M(2n+5)$. Moreover, 
every such $\PSLC$-representation leads to a jumping point.

By Equations~\ref{eqpqr} and \ref{eqH1}, there
are $|n-5|-6$ $\PSLC$-characters of $\Delta(3,5,|n-5|/2)$.
As $M(2n+5)$ has odd degree homology, none of these
are dihedral characters, and they are therefore covered 
twice in $\SLC$. Finally, \cite[Theorem A]{BB}
shows that the jump is two at each of these characters. \qed

\begin{lemma} \label{lem63}
There are $8$ jumping points where $\nu(x)$ is the character
of a reducible representation. These are simple points of
$X$ and each contributes a jump of two.
\end{lemma}

\Pf
If we compare the current case, $n \equiv 5 \pmod{30}$,
with the previous case, $n \not\equiv 5 \pmod{30}$,
we see that there are four fewer irreducible $\PSLC$-characters
of $\Delta(3,5,|n-5|/2)$ ($|n-5|-6$ versus $|n-5|-2$).
On the other hand, we now have $4$ reducible characters corresponding
to representations with image $\Z/15$. The strategy is to
show that these also yield jumping points.

Since $n \equiv 5 \mod 30$, then $|n-5|/2 = 15k$ and $H_1(\dn) \cong 
\Z / 15$. As we have mentioned,
$\Apn$ admits primitive $15$th roots of unity as zeroes and they are simple 
zeroes of $\Apn$.

Let $\xi = e^{2 \pi j i /15}$ be a primitive $15$th root of unity and
let $\rho$ be the reducible $\SLC$ representation of $\pi$ induced 
by $$\rho(\mu) = \left( \begin{array}{cc} 
e^{\pi j i/15} & 0 \\
0 & e^{- \pi j i /15} \end{array} \right). $$
Then $\rho(\mu)^{15} = \pm \I$ and $\rho$ is a cover of one of 
the reducible $\PSLC$ representations of $\dn$ projecting onto 
$\Z / 15$.
In other words, we can think of $\rho$ as a representation of
$\pi$ which factors through $M(2n+5)$.
Corresponding to the eight primitive $15$th roots of unity,  
we have eight $\SLC$ characters.
We will show that the jump at each of these characters is two.

Frohman and Klassen~\cite[Theorem 1.1]{FK} show that such a 
representation $\rho$ is the endpoint of an arc of irreducible
representations. So $\rho \in R_i$, a component of the 
$\SLC$-representation variety containing an irreducible
representation. The corresponding character $x = \chi_{\rho}$
lies on the curve $X_i = t(R_i)$.

Since $\xi = e^{2 \pi j i/15}$ is a primitive 15th root of
unity, $\chi_{\rho}(\mu) = 2 \cos(\pi j/15) \neq \pm 2$.
So $Z_x(f_{\mu}) = 0$, $x$ is a non-trivial character,
and, moreover, $x(\pdM) \neq \{ \pm 2 \}$.
(A character is {\em trivial} if $\chi(\pi) \subset \{ \pm 2 \}$. 
See~\cite[Section 3.2]{P}.)
On the other hand, since $\rho$ factors through $M(2n+5)$, 
$\rho(2n+5) = I$ and $Z_x(f_{2n+5}) > 0$. So
$Z_x(f_{2n+5}) > Z_x(f_{\mu})$ and there is a jump 
at $x$. 

Now, Proposition 1.5.2 of~\cite{CGLS} shows that there is a
non-abelian representation ${\rho}' \in R_i$ with character 
$x$ and ${\rho}'(2n+5) = \pm \I$. Since $x(\pdM) \neq 
\{ \pm 2 \}$, we see that ${\rho}'(\pdM) \not\subset \{ \pm \I \}$.
Finally, as in the proof of Lemma~\ref{lem51}, we can argue that
$H^1(\pi_1(M(2n+5)); {\slC}_{Ad {\rho}'}) = 0$.
These are then simple points of $X$ (\cite[Section 4]{BZ}) 
and each provides a 
jump, $Z_x(f_{2n+5}) - Z_x(f_{\mu})$, of two (\cite[Theorem A]{BB}). \qed

Combining Lemmas~\ref{lem62} and \ref{lem63} we see that we
have also proved Proposition~\ref{prp2n5} in the case 
$n \equiv 5 \pmod{30}$. This completes the proof of the 
proposition. \qed

\section{Proof of Main Theorem}

In this section we prove

\setcounter{section}{1}
\setcounter{theorem}{5}

\begin{theorem} The $\SLC$ character variety of the 
hyperbolic $(-2,3,n)$ pretzel knot $K_n$ contains a curve
of reducible characters and a norm curve $X_0$. If $n \mid 3$,
there is in addition an $r$-curve $X_1$ with $r = 2n+6$ and
$s_1 = 2$. The Culler-Shalen norm $\| \cdot \|_0$ 
for the norm curve is as follows.

If  $3 \nmid n$, then $s_0 = 3(|n-2| -1)$ and
$$ \|\gamma\|_0 = 2 [ \Delta (\gamma,16) + 
2 \Delta( \gamma, \frac{n^2 - n - 5}{ \frac{n-3}{2}}) 
+ \frac{n-5}{2} \Delta( \gamma, 2n+6)]$$
when $n \geq 7$ and 
$$ \|\gamma\|_0 = 2 [  \Delta (\gamma,10)  
+ \frac{1-n}{2} \Delta( \gamma, 2n+6) + 
\Delta( \gamma, 2(n+1)^2 / n) ]$$
when $n \leq -1$.
 
If $3 \mid n$, then 
$s_0 = 3|n-2| -5$ and 
$$ \|\gamma\|_0 = 2 [ \Delta (\gamma,16) + 
2 \Delta( \gamma, \frac{n^2 - n - 5}{ \frac{n-3}{2}}) 
+ \frac{n-7}{2} \Delta( \gamma, 2n+6)]$$
when $n \geq 7$ and 
$$ \|\gamma\|_0 = 2 [  \Delta (\gamma,10)  
- \frac{n+1}{2} \Delta( \gamma, 2n+6) + 
\Delta( \gamma, 2(n+1)^2 / n) ]$$
when $n \leq -1$.  
\end{theorem}

\setcounter{section}{6}
\setcounter{theorem}{0}

\Pf
By Lemma~\ref{lemBZ62}, the 
Culler-Shalen seminorms can be written in terms of 
the boundary slopes $\beta_j$. Using the methods of 
\cite{HO,Du} there are four boundary slopes: $\beta_1 = 0$, $\beta_2 = 2n+6$,
$\beta_3 = 16$ (respectively $10$), and
$\beta_4 = \frac{n^2 -n -5}{\frac{n-3}{2}}$
(respectively $2(n+1)^2/n$) when $n \geq 7$ (respectively
$n \leq -1$). Thus,
$$ \| \gamma \|_i = 2 {[} a^i_1 \Delta( \gamma, \beta_1) + 
a^i_2 \Delta( \gamma, \beta_2) + a^i_3 \Delta( \gamma, \beta_3) + 
a^i_4 \Delta( \gamma, \beta_4) {]},$$
and finding the seminorms comes down to determining the 
non-negative integers $a^i_j$ ($j = 1,2,3,4$). We will 
frequently omit the $i$ super- and subscripts in the following.

Propositions~\ref{prpS}, \ref{prp2n4}, and \ref{prp2n5} imply
the following inequalities for each Culler-Shalen
seminorm $\| \cdot \|$.  
\begin{eqnarray} \label{eqCS}
\| \mu  \|  & = & s \leq 3(|n-2| -1) \mbox{; } \nonumber\\
s \leq \| 2n + 4 \| & \leq & s + 3(|n-6| -1) \mbox{ and }  \\
s \leq \| 2n + 5 \| & \leq & s + 4(|n-5| -2).  \nonumber
\end{eqnarray}

\begin{lemma} \label{lemnorm}
If $n \geq 7$ or $n \leq -11$ there is exactly one norm
curve $X_0$ in the character variety $X$.

If $n \geq 7$, the Culler-Shalen norm on $X_0$ is given by
the coefficients $a_1 = 0$, $a_3 = 1$, $a_4 = 2$, and
$0 \leq a_2 \leq (n-5)/2$.

If $n \leq -11$, the set of  
coefficients is one of the following
three types.
\begin{enumerate}
\item 
$ a_1 = 0 \mbox{, } a_3 = a_4 = 1 \mbox{ and } 0 \leq a_2 \leq (1-n)/2. $
\item
$a_1 = a_4 =1 \mbox{, } a_3 = 0 \mbox{ and } 0 \leq a_2 \leq (1-n)/2.$
\item 
$a_1=1 \mbox{, } a_3 = a_4 = 0 \mbox{ and } 0 < a_2 \leq (n+25)/2.$
\end{enumerate}
Moreover, Type 3 can occur only if $n \geq -23$.
\end{lemma}

\Pf
First, suppose $n \geq 7$. Then Equations~\ref{eqCS} 
become
\begin{eqnarray}
2{[} a_1 + a_2 + a_3 + \frac{n-3}{2}a_4 {]} & = & s \leq 3(n-3) \mbox{; } 
 \label{eqno1p}\\
s \leq 2 {[} (2n+4)a_1 + 2 a_2 + (2n - 12) a_3 + a_4 {]} & \leq & s + 3(n-7) 
\mbox{ and } \label{eqno2p}\\
s \leq 2 {[} (2n+5)a_1 + a_2 + (2n - 11) a_3 + \frac{n-5}{2} a_4 {]}
 & \leq & s + 4(n-7). \label{eqno3p} 
\end{eqnarray}

It will be useful to subtract $s$ from each of the last two equations:
\begin{eqnarray}
0 \leq (2n+3)a_1 + a_2 + (2n - 13) a_3 - \frac{n-5}{2}a_4 & \leq &  
3(n-7)/2,  \label{eqno2pp}\\
0 \leq (2n+3)a_1 + (2n-12)a_3 - a_4 & \leq & 2(n-7) \label{eqno3pp} 
\end{eqnarray}

Since $a_i \geq 0$, Equation~\ref{eqno1p} implies $a_4 \leq 3$.
Moreover, in order to have a norm (rather than a seminorm),
we would need at least two of the $a_i > 0$. This condition
further restricts $a_4 \leq 2$. 

Given $a_4 \leq 2$, Equation~\ref{eqno3pp} implies 
$(2n+3)a_1 \leq 2(n-6)$ so that $a_1 = 0$. Then, the same equation
implies $a_3 \leq 1$.  We will argue that, in fact, $a_4 = 2$ and 
$a_3 = 1$.

Suppose instead that $a_4 \leq 1$. Since $a_1 = 0$, Equation~\ref{eqno3pp}
becomes $(2n-12)a_3 \leq 2n-13$ so that $a_3 = 0$. This is a contradiction
since if $a_1$ and $a_3$ are both zero, then Equation~\ref{eqno3pp}
in fact says $a_4 = 0$ as well, and only $a_2$ is non-zero. This would
mean that $\| \cdot \|$ is not a norm. 

Therefore, $a_4 = 2$. Since $a_1 = 0$, Equation~\ref{eqno3pp}
implies that $a_3 > 0$. Thus $a_3 = 1$. Finally, given these
values, Equation~\ref{eqno2pp} can be rearranged to see that
$0 \leq a_2 \leq (n-5)/2$. This implies $s = 2n - 4 + 2a_2$,
$\|2n+4 \| = s + 2(n-8) + 2a_2$ and $\| 2n+5 \| = s + 4(n-7)$.

Suppose there were two norm curves, $X_1$ and $X_2$. 
Then each would have norm as
described in the previous paragraph. In particular 
$s_1, s_2 \geq 2n-4$. But then $S \geq s_1 + s_2 > 3(|n-2|-1) = 3(n-3)$
which contradicts Proposition~\ref{prpS}. Therefore, there can
be at most one norm curve. On the other hand, since the $(-2,3,n)$ 
pretzel knot is hyperbolic, we know that there is a norm curve in
its character variety, namely the canonical curve.
Therefore, there is exactly one
norm curve when $n \geq 7$. Moreover, the coefficients for that 
curve are $a_1 = 0$, $a_3 = 1$, $a_4 = 2$, and 
$0 \leq a_2 \leq (n-5)/2$. 

For $n \leq -3$ there are several possible solutions to Equations~\ref{eqCS}.
(Although Lemma~\ref{lemnorm} refers to $n \leq -11$, we include
solutions for $n \leq -3$ for future reference. On the other 
hand, solutions for $n = -1$ are not necessary since the $(-2,3,-1)$ pretzel
knot is a twist knot and it's Culler-Shalen seminorms were worked out
in \cite{BMZ}.)

To be precise, there are four possible solutions which lead to
a norm (rather than just a seminorm).
\begin{enumerate}
\item 
$ a_1 = 0 \mbox{, } a_3 = a_4 = 1 \mbox{ and } 0 \leq a_2 \leq (1-n)/2. $
Then
$$s = 2(1-n) + 2a_2 \mbox{,    } 
 \|2n + 4\| = 4(4 -  n) + 4a_2 = s + 2(7 - n) + 2a_2  $$
$$  \mbox{ and }
\|2n + 5 \| = 2(7 - 3n) + 2a_2 = s + 4(3 - n).$$
\item
$a_1 = a_4 =1 \mbox{, } a_3 = 0 \mbox{ and } 0 \leq a_2 \leq (1-n)/2.$
Then $$s = 2(1-n) + 2a_2 \mbox{, } \|2n + 4\| = -4(n + 1) + 4a_2 = s -2(n +3) + 2a_2  $$
$$ \mbox{ and } \|2n + 5 \| = -6(n+1) + 2a_2 = s -4(n + 2).$$

\item If $n \geq -23$,
$a_1=1 \mbox{, } a_3 = a_4 = 0 \mbox{ and } 0 < a_2 \leq (n+25)/2.$
Then
$$s =  2 + 2 a_2 \mbox{, } \|2n+4\| = -2(2n+4) + 4 a_2 = s -2(2n+5) + 2a_2 $$
$$\mbox{ and } \|2n+5 \| = -2(2n+5) + 2a_2 = s -2(2n+6).$$

\item If $n = -3$,
$ a_1 = a_4 = 0 \mbox{ and } a_2 = a_3 =1.$
Then
$$s = 4 \mbox{, } \| 2n +4 \| = 28 \mbox{ and } \| 2n+5 \| = 24. $$
\end{enumerate}

Since $s \geq 2(1-n)$ for curves of Type 1 and 2, and since 
$2 \times 2(1-n) > S$, there is at most one norm curve of Type 1 or
2 in the character variety. To complete the argument for $n \leq -11$,
we must show that there cannot be a norm curve of Type 3 in the 
presence of any other norm curve. 

For example, suppose there were norm curves $X_1$ and $X_2$ of Types 1 and 3 
respectively. Then, 
 $\| 2n+5 \|_1 = s_1 + 4(3-n)$  which implies all the 
jumping points for $2n+5$ surgery lie on the curve $X_1$
(see proof of Proposition~\ref{prp2n5}). 
We've shown that the
jumping points are simple (see Lemmas~\ref{lem62} and \ref{lem63}), 
so they cannot lie on any other
curve. Thus, on $X_2$, $2n+5$ surgery must have minimal norm:
$\| 2n+5 \|_2 = s_2$. This is a contradiction.

Similar arguments show that combinations of Type 3 with Type 2 or
Type 3 with Type 3 are also not feasible. Therefore, there can be 
at most one norm curve when $n \leq -11$. On the other hand, 
since the $(-2,3,n)$ pretzel knot is hyperbolic, there is at least one 
norm curve, namely the canonical curve.
So there is exactly one norm curve when $n \leq -11$
and its norm is of Type 1, 2, or 3. \qed

Lemma~\ref{lemnorm} largely completes the proof of the proposition
in case $n \geq 7$ or $n \leq -11$. What remains is to analyze
the $r$-curves. Recall that $r$ must be an integral boundary 
slope and that the Culler-Shalen seminorm on an $r$-curve is
$\| \gamma \|_i = s_i \Delta(\gamma, r).$ 

\begin{lemma} \label{lemrcrv}
Let $n \geq 7$ or $n \leq -11$. If $3 \mid n$, there is a
unique $r$-curve with $r = 2n+6$. Moreover, the minimal norm
for this curve is $s = 2$. If $3 \nmid n$, there is no 
$r$-curve.
\end{lemma}

\Pf
We will argue that $K_n$ admits an $r$-curve only if $r = 2n+6$.
For this, we'll use the graph manifold structure of $M(2n+6)$.

\begin{claim} \label{clmgrp}
$M(2n+6) = M_1 \cup M_2$, is the union of two Seifert 
fibred manifolds $M_1$ and $M_2$ along a torus. 
$M_1$ is Seifert fibred over $D^2(2,2)$ and
$M_2$ is Seifert fibred over $D^2(3, |n-3|/2)$.
Moreover, regular fibres intersect once on the boundary torus.
\end{claim}

\begin{figure}
\begin{center}
\epsfig{file=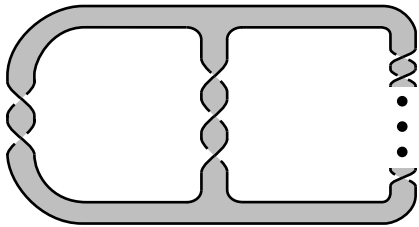}
\end{center}
\caption{\label{fg2n6}%
A spanning surface of $K_n$.}
\end{figure}

\Pf
We can construct $M_1$ by thickening the spanning surface of
Figure~\ref{fg2n6}.
The surface is a punctured Klein bottle so $M_1$ is a twisted 
I-bundle over the Klein bottle. We denote the complementary
manifold $S^3 \setminus N(M_1)$ by $M_2$. We can get a better handle on the
Seifert structure of the $M_i$ ($i=1,2$) using the ideas of Dean~\cite{Dn}.

\begin{figure}
\begin{center}
\epsfig{file=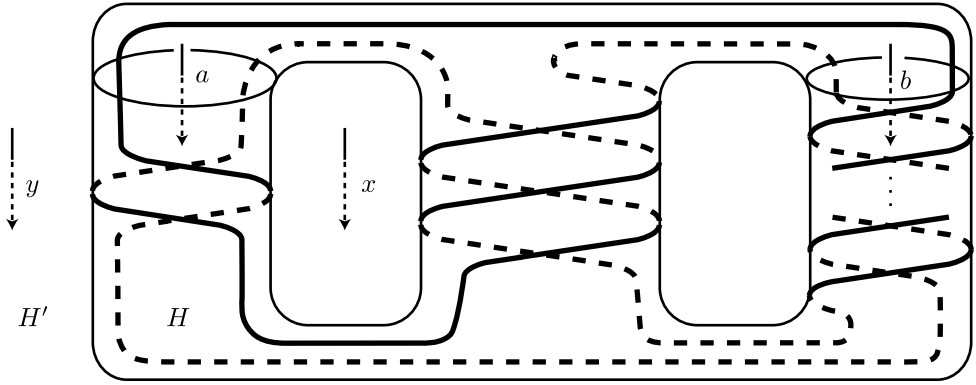}
\end{center}
\caption{\label{fgH}
The handlebody $H$.}
\end{figure}

A key observation~\cite[Lemma 2.2.1]{Dn} is that an irreducible
Haken manifold, like $M_i$, with fundamental group 
$G_{m,n} = \langle x,y \mid x^m y^n \rangle$ 
is Seifert fibred, the base orbifold being a disc with 
cone points of order $m$ and $n$. Since $M_1$ is obtained from the obvious
genus $2$ handlebody $H$ in Figure~\ref{fgH}
by adding a $2$-handle along the knot, we can compute it's
fundamental group. Indeed, with respect to the generators
$a$, $b$ of $\pi_1(H)$, the knot represents the relator 
$b^{-1}ab^{-1}a^{-1}$. Making the change of basis,
$b^{-1}a \to c$, $a \to d^{-1}$, the relator becomes $c^2 d^2$.
Thus $M_1$ is Seifert fibred over $D^2(2,2)$ with $(b^{-1}a)^2$
or $a^2$ representing a regular fibre and fundamental
group $\pi_1(M_1) = \langle c,d \mid c^2d^2 \rangle$.

A similar argument allows us to identify $M_2$ using the generators
$x$ and $y$ of the complementary handlebody $H'$ (see Figure~\ref{fgH}).
In this context, 
the knot represents the word $yxy^{(n-1)/2} x y^{(n-1)/2}x$. After the 
change of basis $y^{(n-1)/2} x \to w$, $y^{-1} \to z$, the word
becomes $z^{(n-3)/2}w^3$. Thus $M_2$ is Seifert fibred over
$D^2(3, |n-3|/2)$ and has fundamental group
$\pi_1(M_2) = \langle w,z \mid z^{(n-3)/2}w^3 \rangle$.
Moreover, a regular fibre corresponds to $(xy^{(n-1)/2})^3$.

\begin{figure}
\begin{center}
\epsfig{file=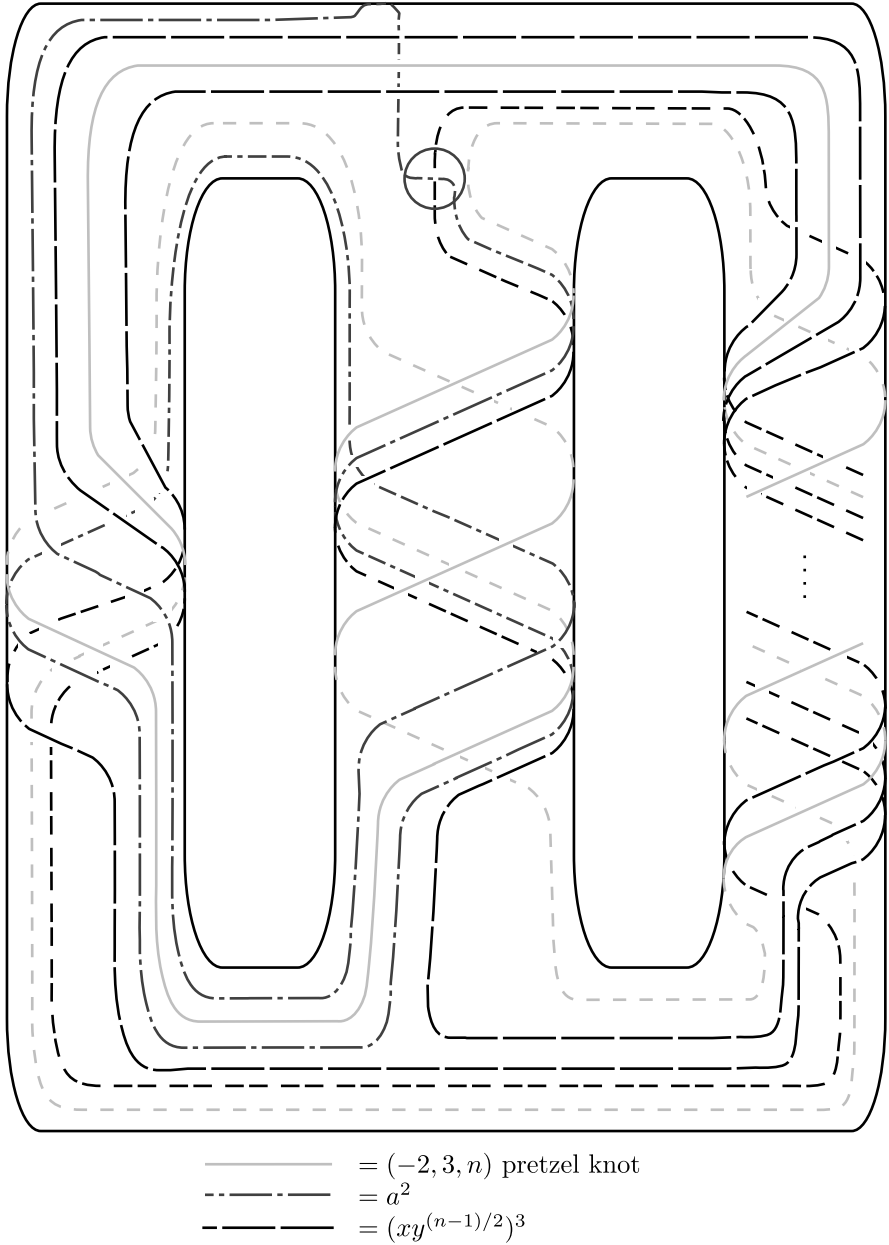}
\end{center}
\caption{\label{fgfbr}
The regular fibres $a^2$ and $(xy^{(n-1)/2})^3$ intersect once
inside the circle.}
\end{figure}

We can now argue that the fibres intersect once on the common
boundary of $M_1$ and $M_2$. Indeed, Figure~\ref{fgfbr} 
shows how a fibre of $M_1$ representing $a^2$ and a fibre of 
$M_2$ representing $(xy^{(n-1)/2})^3$ have intersection 
number one. Note also that the $M_1$ fibre $a^2$ becomes
$y^{-2} x^{-1} = y^{(n-5)/2} (xy^{(n-1)/2})^{-1} = z^{(5-n)/2}w^{-1}$
in $\pi_1(M_2)$ whereas the $M_2$ fibre $(xy^{(n-1)/2})^3$
goes to $b^{-1}ab^{-1}ab^{-1}a = (b^{-1}a)^3 a^{-1} = c^3d$. \qed (Claim)

\begin{claim} \label{clm2n6} If $3 \mid n$,
the $\PSLC$-character variety $\bX(M(2n+6))$ 
contains exactly one curve. If $3 \nmid n$, then 
$\dimC{\bX(M(2n+6))} = 0$.
\end{claim}

\Pf (of Claim)
We will argue that irreducible $\PSLC$-characters of $M(2n+6)$
are either isolated or else factor through $\Z/ 2 \ast \Z / 3$.
The result then follows from \cite[Example 3.2]{BZ2}.

An irreducible representation $\brho : M(2n+6) \to \PSLC$ will
be non-abelian or else have image $\Z /2 \oplus \Z /2$ . On the
other hand, if it's abelian, it also factors through
$H_1(M(2n+6); \Z) \cong \Z / (2n+6)$ and there is no cyclic
group which contains $\Z /2 \oplus \Z /2$. 
Therefore, if $\brho$ is irreducible, it's also non-abelian.
Let $\brho_i : \pi_1(M_i) \to \PSLC$ ($i = 1,2$) be the 
induced representations and let 
$h_i \in \pi_1(M_i)$ be the class of a regular fibre.
A little algebra shows that if $\brho_i$ is non-abelian,
then $\brho_i(h_i) = \pm \I$ (see \cite[Claim 5.2.3]{Ma1}).  

\begin{sub}
If $\brho(\pi_1(T)) \not\subset \{ \pm \I \}$, then
$\chi_{\brho}$ is isolated in $\bX(M(2n+6))$.
\end{sub}

\Pf (of Subclaim)
Let us assume $\brho(\pi_1(T)) \not\subset \{ \pm \I \}$. We wish to 
show that $\chi_{\brho}$ is then isolated in $\bX(M(2n+6))$. 
Since regular fibres intersect once on $T$, their images
generate $\brho(\pi_1(T))$. Therefore, 
in order to satisfy $\brho(\pi_1(T)) \not\subset \{ \pm \I \}$,
at least one $\brho_i$ is abelian with $\brho_i(h_i) \neq \pm \I$.

For example, suppose $\brho_2$ is abelian and $\brho_1$ is not.
As above, $\brho_1(h_1) = \pm \I$. Since the glueing torus $T$ contains
regular fibres, we can assume $h_1 \in \pi_1(T)$. As the $\brho_i$'s
agree on the intersection $\pi_1(T)$, $\brho_2(h_1) = \pm \I$ as well.
However, we've seen earlier (in the proof of Claim~\ref{clmgrp}) 
that a regular fibre $h_1$ represents the
word $z^{(5-n)/2}w^{-1}$ in $\pi_1(M_2)$. Since this word is killed,
$\brho_2$ factors through
\begin{eqnarray*}
\pi_1(M_2) / \langle h_1 \rangle  = 
\langle w,z \mid w^3 z^{(n-3)/2}, z^{(5-n)/2} w^{-1} \rangle
                   = \langle z \mid z^{6-n} \rangle
\end{eqnarray*}
which is cyclic of order $|n-6|$.

This means that $\brho_2(h_2)$ is of finite odd order (remember that 
$n$ is odd, so that $n-6 \neq 0$). We can conjugate so that
$\brho_2(h_2) = \Pmat{\eta}{0}{0}{1/\eta}$ with $\eta \neq \pm 1, \pm i$.
Now, since $\brho_1(h_1) = \pm \I$, $\brho_1$ factors
through the orbifold group $\pi_1^{\mbox{orb}}(\B_1) =
\langle c,d \mid c^2, d^2 \rangle$
where $\B_1$ is the base orbifold of $M_1$. Again, 
$\brho_1(h_2) = \brho_2(h_2)$ is of finite order dividing 
$|n-6|$ and represents the word $c^3d$. Thus
$\brho_1$ factors through $\langle c,d \mid c^2, d^2,
(cd)^{|n-6|}\rangle$ which is dihedral of order $2|n-6|$. Also, 
$\brho_1(h_2) = \brho_2(h_2) = \Pmat{\eta}{0}{0}{1/\eta}$
is in the image of the cyclic subgroup which is therefore
diagonal. 

We have now given a rather specific description of 
$\brho$. Restricted to $\pi_1(M_2)$, it is cyclic of
order dividing $|n-6|$ and diagonal. Restricted to
$\pi_1(M_1)$ it factors through $D_{2|n-6|}$ with 
the cyclic subgroup having image in the diagonal matrices.
Moreover, $\Pmat{\eta}{0}{0}{1/\eta}$ with $\eta \neq \pm 1, \pm i$
is common to the images of $\pi_1(M_1)$ and $\pi_1(M_2)$. There are
only a finite number of characters consistent with such
a representation. Thus, characters of this form are
isolated in the sense that they cannot form a curve.

Similar arguments apply when $\brho_1$ is abelian
and $\brho_2$ is not or when $\brho_1$ and $\brho_2$ are both
abelian. That is, in each of these cases, we find only a finite
number of isolated  characters which, therefore, 
do not form a curve. \qed (Subclaim)

Given the Subclaim, 
the only way to construct a curve in $\bX(M(2n+6))$ is
by making use of representations $\brho$ which kill the
glueing torus $T$ and therefore factor through
\begin{eqnarray*}
\pi_1(M(2n+6)) / \pi_1(T) & = & (\pi_1(M_1) \ast_{\pi_1(T)} \pi_1(M_2)) / 
 \pi_1(T) \\
 & = & \pi_1(M_1) / \pi_1(T) \ast \pi_1(M_2) / \pi_1(T) \\
 & = & \Z / 2 \ast \Z / g .
\end{eqnarray*}
where $g = \mbox{gcd}(3, |n-3|/2)$. If $3 \nmid n$, then 
$g = 1$ and this is an abelian
representation, contradicting an earlier assumption. So there
is no curve in $\bX(M(2n+6)$ when $3 \nmid n$.
If $g = 3$ (i.e., if $3 \mid n$),
we see that we are looking
at representations of $\Z / 2 \ast \Z / 3$. Since $\bX( \Z / 2 \ast \Z / 3)$
contains exactly one curve (see \cite[Example 3.2]{BZ2}), 
we conclude that this is also the case for $\bX(M(2n+6))$. \qed (Claim)

So, if $3 \nmid n$, $\dimC{\bX(2n+6)} = 0$ and there can be no
$r$-curve with $r = 2n+6$ (compare \cite[Example 5.10]{BZ2}).

On the other hand, 
if $3 \mid n$, there is a unique curve in $\bX(M(2n+6))$.
Moreover, since 
the representation $\brho_{1/2}$ (see \cite[Example 3.2]{BZ2}
or Equation~\ref{eqex32} below) is dihedral, 
this curve contains the character of a dihedral representation.
It follows that the curve is covered by a unique curve
in the $\SLC$-character variety $X(M(2n+6))$ (see \cite[Lemma 5.5]{BZ}). 
Thus there is exactly one $r$-curve, call it $X_1$,
with $r = 2n+6$ when $3 \mid n$.

\begin{claim}
The minimal norm $s_1$ of $X_1$ is two.
\end{claim}

\begin{figure}
\begin{center}
\epsfig{file=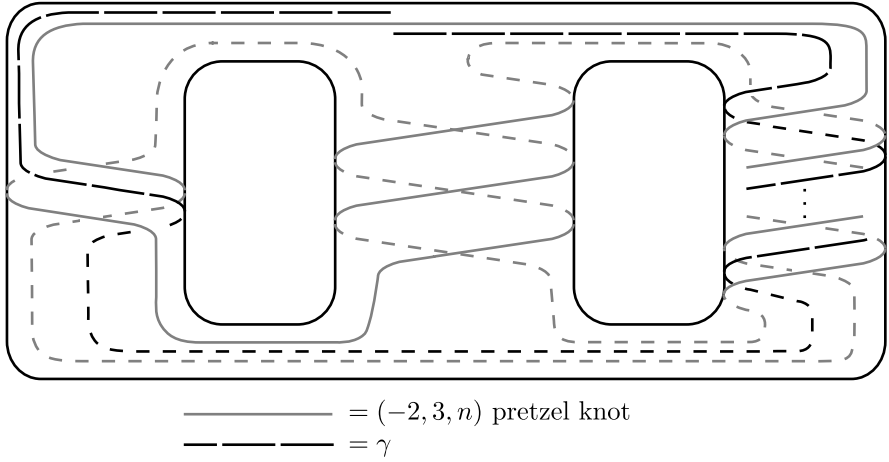}
\end{center}
\caption{\label{fgimmu}
The meridian $\mu$ represented by a curve $\gamma$.}
\end{figure}

\Pf (Claim)
Recall that $s_1 = \| \mu \|_1$ is the degree of $f_{\mu}$ 
(Theorem~\ref{Thcyc}).
So we need to understand the image of $\mu$
under the composition $ \pn \rightarrow \pi_1(M(2n+6)) \rightarrow
\Z/2 \ast \Z/3 = \langle c,d \mid  c^2, d^3 \rangle$. 
We can construct $\mu$ in terms of
a curve $\gamma$ in the genus two surface which connects points on opposite
sides of the knot; see Figure~\ref{fgimmu}.
The idea is that we can break up the meridian as the sum of a loop in
$M_1$ and a loop in $M_2$. We will then show that those project to $c$
and $d$ respectively. The first loop is $\gamma$ plus a small arc joining
the two endpoints of $\gamma$ in the interior of $H$, the obvious 
genus two handlebody. The second loop is $\gamma$ plus a small arc joining
the two endpoints of $\gamma$ and passing through the complementary
genus two handlebody $H'$. 

In $\pi_1(M_1)$, $\gamma$ represents $ab^{-1}$ which is conjugate to 
$b^{-1}a$ and therefore projects onto the generator of 
$\Z / 2  = \pi_1(M_1) / \pi_1(T)$. In $\pi_1(M_2)$, $\gamma$ represents
$xy^{(n-1)/2}$ which projects to the generator of $\Z / 3$. Thus $\mu$
is mapped to $cd$ in $\Z/2 \ast \Z/3$.

Let $\bX_1$ be the unique curve in $\bX(\Z/2 \ast \Z/3)$.
In \cite[Example 3.2]{BZ2}, the authors construct a double 
covering $\C \to \bX_1$ given by mapping $z \in \C$ to the 
character of $\brho_z$:
\begin{equation} \label{eqex32}
\brho_z(c) = \pm \left( \begin{array}{cc} i & 0 \\
                                      0 & -i \end{array} \right)
\mbox{, } 
   \brho_z(d) = \pm \left( \begin{array}{cc} z & 1  \\
                               z(1-z)-1 & 1-z \end{array} \right).
\end{equation}
Since 
$\tr{\brho_z(cd)} = i(2z-1)$, we see that 
$f_{cd}(\chi_{\brho_z}) = -(2z-1)^2 - 4$ has degree $2$.
As this is a double covering of $\bX_1$ by $\C$, the corresponding
character on $\bX_1$ has degree $1$. Finally, lifting to the curve
$X_1$ in $X(M)$ which double covers $\bX_1 \subset \bX(M(2n+6)) 
\subset \bX(M)$, we deduce
$s_1 = \deg{f_{\mu}} = 2$. \qed (Claim)

Thus far, we've constructed an $r$-curve $X_1$ for $r=2n+6$ when $3 \mid n$
and shown that there is none when $3 \nmid n$. It remains to argue
that there are no $r$-curves for the other integral boundary slopes.

For $n \geq 7$, the integral boundary slopes are $0$, $16$, and $2n+6$.
We will show that $0$ and $16$ do not admit $r$-curves.
Recall (Lemma~\ref{lemnorm}) that $\| 2n+5 \|_0 = s_0 + 4(|n-5| -2)$,
i.e., all the jumping points for the $2n+5$ surgery are
on the norm curve $X_0$. So we would have $\| 2n+5 \|_i = s_i$ on any 
$r$-curve $X_i$. However, if $r = 0$ for example, 
$\|2n+5\|_i = s_i \Delta(2n + 5, 0) = (2n+5)s_i$. So there can be
no $r$-curve for $r = 0$. Similarly, there can be no $r = 16$ curve.
Analogous arguments show that when $n \leq -11$,
there is no $r$ curve with $r = 0$ or $10$.
Thus $r=2n+6$ is the only candidate for an $r$-curve
amongst the integral boundary slopes when $n \leq -11$ as well.

This completes the proof of Lemma~\ref{lemrcrv}. \qed 

Thus when $n \geq 7$ or $n \leq -11$, there is exactly one norm curve.
There will be one $r$-curve when $3 \mid n$ and otherwise there are no 
additional curves containing irreducible characters.
Since the set of reducible characters forms a complex line, we
see that $X(K_n)$, the character variety of 
the knot $K_n$, consists of two (three) curves when $3 \nmid n$ 
($3 \mid n$) and $n \geq 7$ or $n \leq -11$. 
To complete the proof of the Main Theorem for these $n$, 
we need only verify that the Culler-Shalen seminorms are as stated.

If $n \geq 7$ we know the norm on the norm curve $X_0$ up to the 
coefficient $a_2$. If $3 \nmid n$, then there is no $r$-curve and
the norm curve $X_0$ is the only one contributing to the
total norm $\| \cdot \|$. In particular $s_0 = S = 3(|n-2|-1) = 3(n-3)$.
Therefore, $a_2 = (n-5)/2$ and the Culler-Shalen norm on $X_0$ is
as stated in the theorem. If $3 \mid n$, there is also an $r$-curve
$X_1$ with $s_1 = 2$. In this case $S = s_0 + s_1 \Rightarrow
s_0 = S - s_1 = 3(n-3) -2$. This implies $a_2 = (n-7)/2$.

If $n \leq -11$, we have several candidates (Types 1, 2, and 3 of 
Lemma~\ref{lemnorm}) for the norm $\|\cdot \|_0$ on the norm curve $X_0$.
If $3 \nmid n$, there are no other curves and $\| \cdot \|_0$
{\em is} the total norm. In particular, it must satisfy 
Propositions~\ref{prpS}, \ref{prp2n4}, and \ref{prp2n5}. The 
only possibility is $\| \cdot \|_0$ of Type 1 with $a_2 = (1-n)/2$.
If $3 \mid n$, the total norm is a combination of 
$\| \gamma \|_1 = 2 \Delta(2n+6, \gamma)$ on $X_1$ and the norm
$\| \cdot \|_0$ on the norm curve $X_0$. Again, the only 
choice satisfying 
Propositions~\ref{prpS}, \ref{prp2n4} and \ref{prp2n5} is
$\| \cdot \|_0$ of Type 1 with $a_2 = -(n+1)/2$. 

This completes the proof of the Main Theorem for the cases 
$n \geq 7$ and $n \leq -11$. There remain the cases $n = -9, -7, \ldots -1$. 

\vspace{12 pt}

\noindent${\mathbf {n = -9}}$
Here $3 \mid n$, so, following the reasoning of Lemma~\ref{lemrcrv},
there is an $r$-curve $X_1$ with $r = 2n + 6 = -12$ and 
$s_1 = 2$. Moreover, there is no $r$-curve with
$r = 0$ or $10$.

As for the norm curves, if there's a norm curve of 
Type 1 or 2, it's the only norm curve. For example, a norm curve $X_0$
of Type 1 has $\| 2n+5 \|_0 = s_0 + 4(3-n)$ so that all the jumping
points for this surgery are on $X_0$. This means any other norm 
curve would have minimal norm for this surgery slope: $\| 2n+5 \|_i = s_i$.
However, this is not true for any norm curve, be it of Type 1, 2, or 3..

On the other hand, we need a new type of argument to show that
two (or more) curves of Type 3 is also not a feasible arrangement. 
Suppose then 
that there were two Type 3 norm curves $X_0$ and $X_2$. 
Since $S = 30$ and $s_1 = 2$, we see that 
$$28 \leq s_{0} + s_{2} = 4 + 2(a_{2}^{0} + a_{2}^{2}),$$
where we have given the $a_2$'s 
superscripts showing which curve they come from.
This implies $a_{2}^{0} + a_{2}^{2} \geq 12$. 
But then  $$ \| 2n+4 \|_{0} + \|2n+4 \|_{2} = 
s_{0} + s_{2} + -4(2n+5) + 2(a_{2}^{0} + a_{2}^2) \geq
s_{0} + s_{2} + 52 + 24$$
which contradicts the equation $\| 2n+4 \|_T  =  S + 3(|n-6| -1) = S + 42$.
Thus, we see that there is exactly one norm curve and one $r$-curve when
$n = -9$. Since these must combine to give the total norm, 
the norm curve is of Type 1 with $a_2 = 4 = -(n+1)/2$.

\vspace{12 pt}

\noindent${\mathbf {n = -7}}$
In this case there is no $r$-curve for $r = 2n + 6 = -8$. 
By examining the norm of the $2n+5 = -9$ slope, we see that 
if there is a norm curve of Type 1, it is the only curve 
in $X$ containing an irreducible character (in particular,
there are no $r$-curves). 

Similarly, if there is a norm curve of Type 2, then
there is no $r$-curve with $r = 10$. 

However, we cannot 
immediately eliminate the possibility of an $r$-curve for $r = 0$.
Indeed, suppose that there were a Type 2 norm curve $X_0$,
together with an $r$-curve $X_1$ with $r=0$.
For the norm curve $X_0$,
$$s_0 = 2(1-n) + 2a_2 = 16 + 2a_2 \mbox{, } 
\| 2n+4 \|_0 = \|-10\|_0 = s_0 -2(n + 3) + 2a_2 = s_0 + 8 + 2a_2 $$
$$ \mbox{ and } \| 2n+5 \|_0 = \|-9\|_0 = s_0 -4(n+2) = s_0 + 20.$$
On the other hand, for the $r$-curve $X_1$,
$\|2n+4 \|_1  = \|-10 \|_1 \leq s_1 + 3(5-n)  - 8 = s_1 +28$ and also
$\|2n+4 \|_1 = \| -10 \|_1 = s_1 \Delta(-10,0) = 10 s_1$.
This implies 
\begin{eqnarray*}
10 s_1 & \leq & s_1 + 28 \\
\Rightarrow 9 s_1 & \leq & 28 \\
\Rightarrow s_1 & \leq &  28/9 
\end{eqnarray*}

Since $s_1$ is an even integer, we see that $s_1 = 2$. Similarly, an
examination of the $-9$ slope also leads us to the conclusion 
that $s_1 = 2$. So we cannot eliminate the possibility of an $r$-curve with
$r = 0$ directly as we did earlier. We need a new type of argument
to handle this situation. We need to analyze the possible combinations
of curves.

For example, suppose $X$ contained a Type 2 norm curve $X_0$
and one $r$-curve $X_1$
for $r=0$ and {\em no other} norm or $r$-curves. Then $\|-9\|_0 = s_0 + 20$
and 
\begin{eqnarray*}
\|-9\|_1 & = & s_1 \Delta(-9,0) \\
         & = & 9 s_1  \\
         & = & s_1 + 8 s_1 \\
         & = & s_1 + 16.
\end{eqnarray*}
So 
\begin{eqnarray*}
\| -9 \|_T & = & \| -9 \|_0 + \|-9 \|_1 \\
           & = & s_0 + s_1 + 36 \\
           & = & S + 36 \\
           & < & S + 40 \\
           & = & S + 4(|n-5| -2). 
\end{eqnarray*}

Thus, if we assume that these are the only two curves, we see that 
we cannot account for all the jumping points associated with the 
$-9$ slope. Therefore, this is not a possible configuration for 
$X$. By analyzing the possible combinations of norm curves
and $r$-curves in this way, we see that the only possibility consistent
with Propositions~\ref{prpS}, \ref{prp2n4}, and \ref{prp2n5} is 
that there is exactly one norm curve of Type 1 with
$a_2 = 4 = (1-n)/2$ and no $r$-curves.

\vspace{12 pt}

\noindent${\mathbf {n = -5}}$
A similar analysis shows that $X$ contains no
$r$-curves and exactly
one norm curve of Type 1 with $a_2 = 3 = (1-n)/2$.

\vspace{12 pt}

\noindent${\mathbf {n = -3}}$
Since $3 \mid n$, 
we know (see Lemma~\ref{lemrcrv})
that there is an $r$-curve $X_1$ for $r = 2n + 6 = 0$ with
$$ s_1 = 2 \mbox{, } \| 2n + 4 \|_1 = \| -2 \|_1 = s_1 + 2 \mbox{ and }
\|2n + 5 \|_1 = \| -1 \|_1 = s_1.$$
If we follow the same strategy as in the previous cases
we find that there are two possible configurations:
\renewcommand{\theenumi}{\Roman{enumi}}
\begin{enumerate}
\item
In addition to the $r$-curve there is one Type 1 norm curve $X_0$ with
$$s_0 = 10 \mbox{, } \| -2 \|_0 = s_0 + 22 \mbox{ and } 
\| -1 \|_0 = s_0 + 24.$$
\item
Here there is a Type 2 norm curve $X_0$:
$$s_0 = 8 \mbox{, } \| -2 \|_0 = s_0  \mbox{ and } 
\| -1 \|_0 = s_0 + 4$$
as well as an additional $r$-curve $X_2$ with $r = 10$:
$$ s_{2} = 2 \mbox{, } \| -2 \|_{2} = s_{2} + 22 \mbox{ and }
\| -1 \|_{2} = s_{2} + 20.$$
\end{enumerate}

Both configurations are consistent with 
Propositions~\ref{prpS}, \ref{prp2n4}, and \ref{prp2n5}: 
$S = 3(|n-2| - 1) = 12$,
$\|-2 \|_T = S + 3(|n-6| -1) = S + 24$, and 
$\| -1 \|_T = S + 4(|n-5|-2) = S + 24$.

In order to show that the second configuration does not arise, 
note that $\bX_2$, the $\PSLC$ analogue of $X_2$, 
would include into $\bX(M(10))$ (see \cite[Example 5.10]{BZ2}).
Now, $M(-1)$ is Seifert fibred over $S^2(3,4,5)$ and the
jumping points for $\|-1\|$ come from the six 
irreducible $\PSLC$-characters of $\Delta(3,4,5)$. If the second
configuration is valid, five of these characters 
are on $\bX_2$ and therefore 
come from representations lying in $\bar{R}(M(10))$. 
We will argue that at least two of them do not.

Indeed two of the characters correspond to representations which factor 
through $\Delta(2,3,5)$ which has order $60$. On the other hand,
if such a representation $\brho$ is also in $\bar{R}(M(10))$, then 
it annihilates both the $10$ and the $-1$ slopes. In other words, 
the kernel of $\brho$ contains an index eleven subgroup of $\pi_1(\pM)$.
Therefore, $\brho(\pi_1(\pM))$ is either $\Z/11$ 
or else trivial. On the other hand, $\brho(\pi_1(\pM))$
also factors through $\Delta(2,3,5)$. Thus $\brho(\pi_1(\pM))$ is
trivial and since $\pi_1(M)$ is normally generated by the peripheral
group, $\brho(\pi_1(M)) = \{ \pm I \}$ as well. This contradicts the
fact that $\brho$ is an irreducible representation. Therefore, 
the irreducible representations which factor through $\Delta(2,3,5)$
are not in $\bar{R}(M(10))$. This shows that the second configuration is
not possible. 

Therefore, we have the first configuration. There is an $r$-curve
and a norm curve of Type 1 with $a_2 =1 =  -(n+1)/2$. 

\vspace{12 pt}

\noindent${\mathbf {n = -1}}$
This knot was treated using different methods in \cite{BMZ}
(where it is identified as the twist knot $K_2$). We saw that there
is one norm curve in the character variety and no $r$-curves.
Moreover, the norm curve corresponds to Type 1 with $a_2 = 1 = (1-n)/2$.

This completes the proof of the Main Theorem. \qed

\section{Applications}

With the Main Theorem in hand, we look at applications. 
Since surgeries which result in a manifold having cyclic
or finite fundamental group are of small norm, we can use
our knowledge of the Culler-Shalen seminorms to understand 
which surgeries might be finite or cyclic. A nice way 
to visualize this connection is via the fundamental polygon.

Given a Culler-Shalen seminorm $\| \cdot \|$ arising from
a norm curve in the character variety of a knot $K$, we call $B$, the disc of 
radius $s$ in $V = H_1(\pM; \R)$, a {\em fundamental polygon} for $K$.
A fundamental polygon is a compact, convex, and 
finite-sided polygon with vertices which are rational multiples of
boundary slopes in $L = H_1(\pM ; \Z).$ It is symmetric ($-B = B$)
and centred at $(0,0)$. 
As we have shown in the Main Theorem, $K_n$ has only one
norm curve $X_0$ which is therefore the canonical curve.
We will refer to the associated polygon as
{\em the} fundamental polygon $B$ of $K_n$. 

By Theorems~\ref{Thcyc} and \ref{Thfin}, a cyclic or finite 
surgery slope either has norm bounded by $\mbox{max}(2s_i, s_i +8)$ or
else is a boundary slope.
On the other hand, for a small knot like $K_n$, a cyclic or
finite surgery cannot occur on a boundary slope.

\begin{lemma} \label{lefinbdy}
If $M$ is small and $\al$ is a boundary slope, then $M(\al)$ is not
cyclic or finite.
\end{lemma}

\Pf By \cite[Theorem 2.0.3]{CGLS}, $M(\al)$ is not finite, and
it is cyclic only if $M(\al) \cong S^1 \times S^2$.
However, Gabai~\cite{Gb} has shown that, amongst knots in $S^3$, 
only slope $0$ surgery on the trivial knot can produce $S^1 \times
S^2$. \qed

Therefore, as long as $s_0 > 8$, all cyclic and finite
surgeries of $K_n$ will lie within $2B$, the norm disk of
radius $2s_0$.

\setcounter{section}{1}
\setcounter{theorem}{6}

\begin{prop}
If the $(-2,3,n)$ pretzel knot $K_n$ admits a non-trivial
cyclic or finite surgery, then one of the following holds.
\begin{itemize}
\item $K_n$ is torus, in which case $n = 1$,$3$, or $5$,
\item $n =7$, in which case $18$ and $19$ are cyclic
fillings while $17$ is a finite, non-cyclic filling, or
\item $n=9$, in which case $22$ and $23$ are finite, non-cyclic
fillings.
\end{itemize}
\end{prop}

\setcounter{section}{7}
\setcounter{theorem}{1}

\Rmk
A torus knot admits an infinite number of cyclic fillings.
A cyclic filling of a non-trivial knot in $S^3$ is necessarily finite
\cite{Gb}.
The finite and cyclic surgeries of $K_7$ and $K_9$ were discovered
by Fintushel and Stern (see~\cite{FS}) and Bleiler and Hodgson~\cite{BH}.
The content here is that these are the only non-trivial
finite or cyclic surgeries on the non-torus
members of this family of knots.

\medskip

\begin{figure}
\begin{center}
\epsfig{file=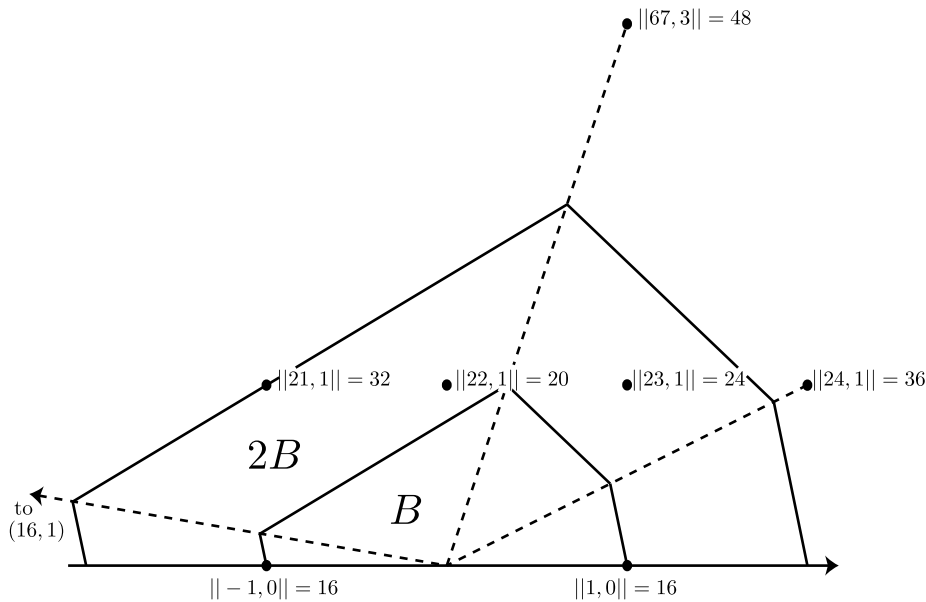}
\end{center}
\caption{The Fundamental Polygon of $K_9$ \label{fg239}}
\end{figure}

\Pf
Let $n \geq 7$. The case $n=7$ is the subject of 
\cite[Example 10.1]{BZ} where it is shown that there are 
exactly three non-trivial finite surgeries. 
For $n=9$, $\max (2s_0, s_0 + 8) = 2s_0$ (Theorem~\ref{thmM})
so any finite or
cyclic surgery slopes must lie in $2B$. However, as we see in 
Figure~\ref{fg239},
the only slopes inside $2B$ are
$21$, $22$, $23$ and $\mu = 1/0$. If $21$ were a finite surgery,
then it would be one with norm $2s_0$ and therefore dihedral
(\cite[Theorem 2.3]{BZ}). However, it cannot be dihedral since
$21$ is odd. Therefore,   
the $(-2,3,9)$ pretzel knot admits exactly two non-trivial
finite surgeries: $22$ and $23$. 

As $n$ increases, the fundamental polygon
for the norm curve maintains the same basic shape but becomes
smaller (the exact coordinates of the polygon are determined by 
the Culler-Shalen norm $\| \cdot \|_0$ of the Main Theorem).
For $11 \leq n \leq 19$, the only slopes inside $2B$ are 
$2n+4$, $2n+5$ and $\mu$ and once $n \geq 21$, only $2n+4$ and $\mu$
remain. However, since $M(2n+4)$ and $M(2n+5)$ are Seifert fibred 
over a hyperbolic orbifold when $n \geq 11$, these are not finite
surgeries. Thus the $(-2,3,n)$ pretzel knots admit no non-trivial 
cyclic or finite surgeries when $n \geq 11$.

\begin{figure}
\begin{center}
\epsfig{file=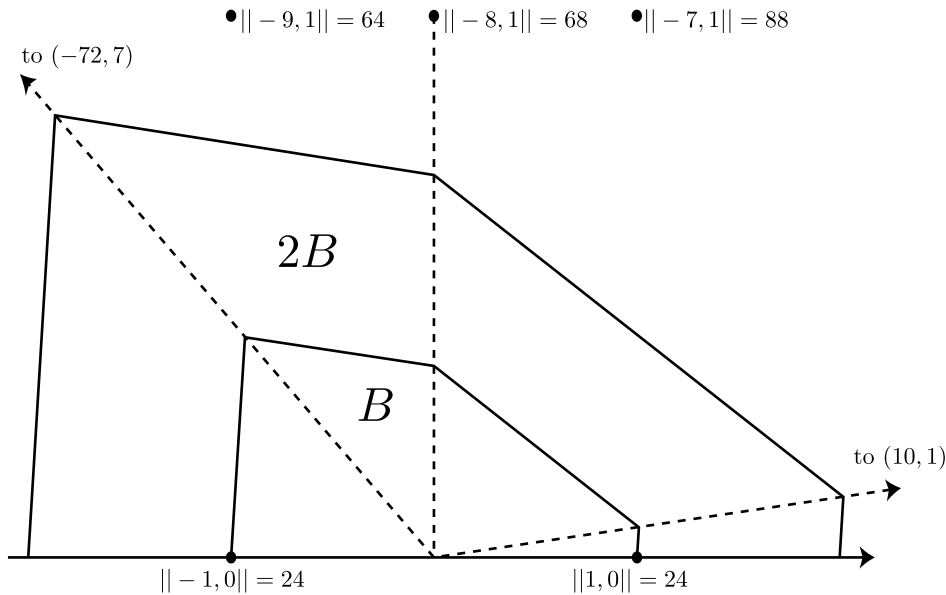}
\end{center}
\caption{The Fundamental Polygon of $K_{-7}$ \label{fg23n}}
\end{figure}

Figure~\ref{fg23n} 
gives the fundamental polygon of the $(-2,3,-7)$ pretzel knot and
illustrates the situation for $n \leq -1$. 
For these knots, $2B$ lies below the line $y = 1$. Thus the
only surgery slope within $2B$ is $\mu$. Since
$s_0 > 8$ for $n \leq -3$, there can be no non-trivial 
finite or cyclic surgeries when $n \leq -3$. For 
$n =-1$, $s_0 = 6$ and any finite or cyclic surgery would have
to lie in the polygon $\frac{14}{6}B$ (i.e., 
$\frac{14}{6} = (s_0 + 8)/s_0$). However, as this polygon
also lies below $y=1$, we again conclude that there are
no non-trivial finite or cyclic surgeries when $n = -1$.

This completes the proof of the proposition. \qed

Note that, as the five finite fillings of the $(-2,3,7)$ and 
$(-2,3,9)$ pretzel knots are not simply-connected, 
this constitutes a proof of Property P for the $(-2,3,n)$ 
pretzel knots. However, these knots are strongly invertible, so 
Property P was already known \cite{BS}.

As a final application, we derive the Newton polygon 
of the $A$-polynomial~\cite{CCGLS} of the $(-2,3,n)$ knots.
Recall that the Newton polygon of a polynomial 
$A = \sum_{(i,j)} b_{i,j} \gl^i \m^j \in
 \Z {[} \gl, \m {]}$
is the convex hull in $\R^2$ of $\{ (i,j) | b_{i,j} \neq 0 \}$.
Boyer and Zhang have shown that 
the Newton polygon $N$ and the fundamental polygon $B$ 
are dual in the following sense. 

\begin{theorem}[Theorem 1.4 of \cite{BZ4}]
The line through any pair of antipodal vertices of $B$ is 
parallel to a side of $N$. Conversely, the line through any
pair of antipodal vertices of $N$ is parallel to a side of $B$.
\end{theorem}

Thus given the fundamental polygon, one can deduce $N$,
at least up to scaling and translation. 
The conventions we use 
are that $N$ meets the $\gl$ and $\m$ axes but lies 
in the first quadrant. The scale is provided by Shanahan's 
width function.

\begin{defn}[Definition 1.2 of \cite{S}] The $p/q$ width
$w(p/q)$
of $N$ is one less than the number of lines of slope $p/q$
which intersect $N$ and contain a point of the 
integer lattice.
\end{defn}

We require that Shanahan's width be half the  
Culler-Shalen norm: $\| q/p \|_0  = 2 w(p/q)$.
With these conventions, the Newton polygon is completely
determined by the norm $\| \cdot \|_0$ of our Main Theorem.

\begin{figure}
\begin{center}
\epsfig{file=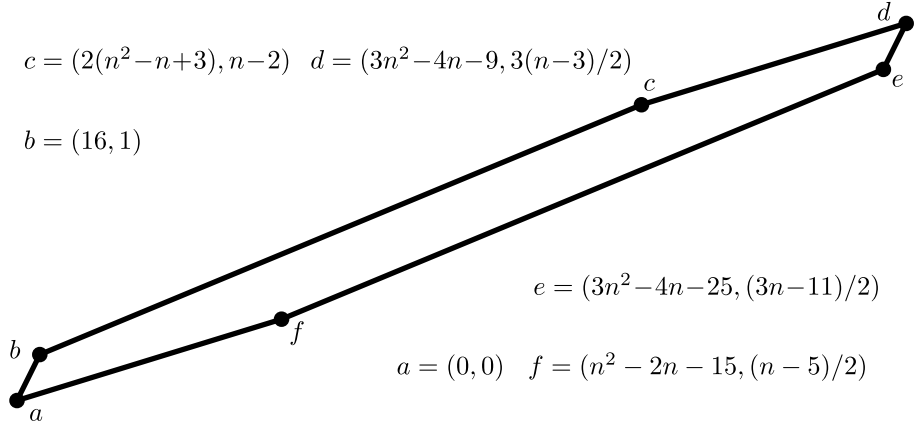}
\end{center}
\caption{The Newton Polygon of $K_n$ ($n \geq 7$ and $3 \nmid n$)
 \label{fgNpos}}
\end{figure}

\begin{figure}
\begin{center}
\epsfig{file=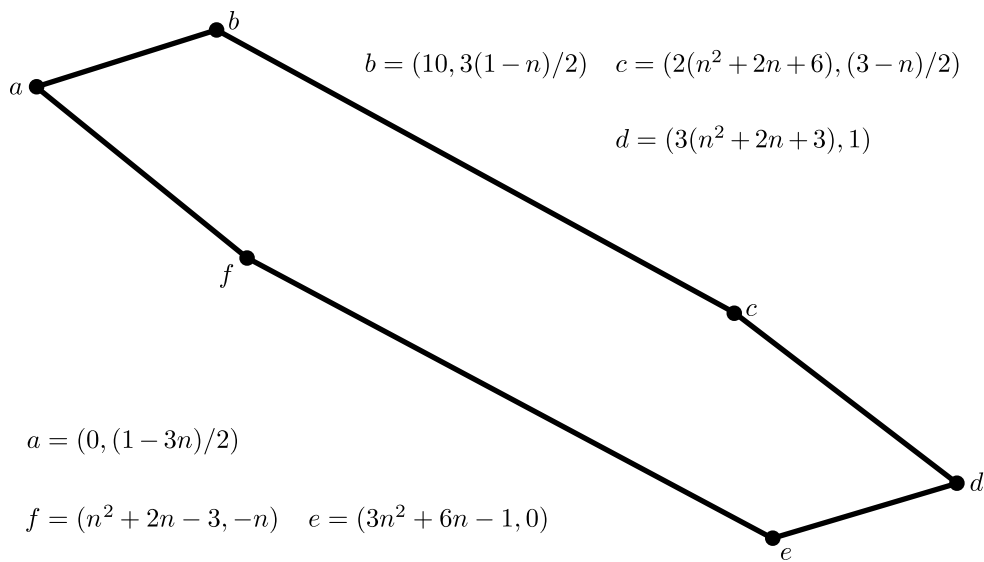}
\end{center}
\caption{The Newton Polygon of $K_n$ ($n \leq -5$ and $3 \nmid n$)
 \label{fgNneg}}
\end{figure}

The vertices of the Newton polygon are 
$$(0, 0), (16,1), (n^2 -2n -15, (n-5)/2), (2(n^2 - n +3), n-2),$$
$$ (3n^2 -4n - 25, (3n-11)/2), (3n^2 -4n-9, 3(n-3)/2)$$
when $n \geq 7$ and $3 \nmid n$ (see Figure~\ref{fgNpos} which may
be compared with \cite[Figure 6]{S}); 
$$(0,0), (16,1), (12(n-7), (n-7)/2), (3(n^2 -6n + 23), n-2),$$
$$ (3n^2 -6n -31,(3n-13)/2), (3(n^2 -2n - 5), (3n-11)/2)$$
when $n = 3k$,  $k \geq 3$;
$$(0, (1-3n)/2), (10, 3(1-n)/2), (n^2 + 2n - 3, -n), $$
$$(2(n^2 +2n + 6), (3-n)/2), (3n^2 +6n -1,0), (3(n^2 + 2n+3),1)$$
when $n \leq -5$ and $3 \nmid n$ (see Figure~\ref{fgNneg});
$$(0, -(3n+1)/2), (10, (1-3n)/2), (n^2 + 4n +3,-n),$$
$$(2(n^2 +2n +6), (1-n)/2), (3n^2 + 8n + 5,0), (3n^2 +8n + 15,1)$$
when $n = 3k$, $k \leq -1$; and
$$(0,0), (0,1), (4,2), (10,1), (14,2), (14,3) \mbox{ when } n = -1.$$

\section*{Acknowledgements}

This work forms part of my Ph.D.\ thesis and I would like 
to thank my supervisor Steven Boyer for his substantial
contributions and indispensable advice. Much of this was written during
a visit to Nihon University in Tokyo and I am grateful to Kimihiko Motegi 
and the Department of Mathematics for their hospitality during my stay.
I am indebted to Kurt Foster 
and Dave Rusin who provided suggestions about cyclotomic roots of the 
Alexander polynomial. 

I am especially grateful to the referee
for a close reading of an earlier version of this paper and
for many concrete 
suggestions which have significantly improved the exposition.

\end{document}